\documentclass[ 11pt] {amsart}	
\usepackage{amssymb}
\usepackage{amscd}
\usepackage{eucal}
\usepackage[dvips]{graphicx}

\newtheorem{theo}{Theorem}[section]
\newtheorem*{theoa}{Theorem A }
\newtheorem*{theob}{Theorem B }
\newtheorem{prop}[theo]{Proposition}

\newtheorem{la}[theo]{Lemma }

\newtheorem{Theorem}{Theorem}
\newtheorem{Corollary}{Corollary}

\theoremstyle{remark}
\newtheorem*{rem}{Remark}

\theoremstyle{definition}

\newcommand{\N}{\mathbb N}

\newcommand{\Q}{\mathbb Q}
\newcommand{\Z}{\mathbb Z}

\newcommand{\A}{\mathbf A}
\newcommand{\F}{\mathbb F}

\newcommand{\we}{\omega(e)}

\newcommand{\tc}{\tilde C}
\newcommand{\mc}{\tilde{\mathcal C}}
\newcommand{\cc}{\mathcal C}
\newcommand{\mi}{\mathcal I}
\newcommand{\tL}{\tilde {\mathbf L}}

\newcommand{\LL} {\mathbf L}
\newcommand{\K} {\mathbf K}
\newcommand{\bc} {\mathbf C^{\bullet}}
\newcommand{\fQ} {\mathbf Q}
\newcommand{\fS} {\mathbf S}
\newcommand{\bp} {\mathbf P}
\newcommand{\hl} {\mathcal L^{\bullet} }
\newcommand{\al} {\mathcal L^{'\bullet}}
\newcommand{\bl} {\mathcal L^{''  \bullet}}
\newcommand{\fu} {\mathbf U}

\newcommand{\rep}{{\mathcal R}}
\newcommand{\AAA}{{\mathcal A}}
\newcommand{\iso}{{\stackrel{\sim}{\rightarrow}}}

\DeclareMathOperator{\Hom}{\rm Hom}

\DeclareMathOperator{\Ind}{\rm Ind}
\DeclareMathOperator{\supp}{\rm supp}
\DeclareMathOperator{\ord}{\rm ord}
\DeclareMathOperator{\Gal}{\rm Gal}

\numberwithin{equation}{section}
\begin{document}
\title{Group cohomology of the universal ordinary distribution}
\author{Yi Ouyang}
\address{
School of Mathematics, University of Minnesota, Minneapolis, 
MN 55455, USA} 
\email{youyang@math.umn.edu}
\begin{abstract}
For any odd squarefree integer $r$, we get complete description of 
the $G_r=\Gal(\Q(\mu_r)/\Q)$ group cohomology of  the universal ordinary 
distribution $U_r$ in this paper.  Moreover,   if M is a fixed integer which divides $\ell-1$
for all prime factors $\ell$ of $r$, we study the cohomology group $H^{\ast}(G_r, U_r/MU_r)$.
In particular,  we explain the  mysterious construction of  the elements $\kappa_{r'}$ for 
$r'|r$  in Rubin\cite{Rubin1},  which come exactly from a certain $\Z/M\Z$-basis 
of the cohomology group $H^{0}(G_r, U_r/MU_r)$ through an evaluation map. 
\end{abstract}
\maketitle
\section{Introduction}

Let $\{[a]: a\in \Q/\Z$ \} be a basis for a free 
abelian group $\A$. Then the (dimension $1$) universal ordinary distribution 
$U_r$ of level $r$ for any positive integer $r$ is given by 
\[ U_r=\frac{<[a]: a\in \frac{1}{r}\Z/\Z>}{<[a]-{\displaystyle \sum_{\ell b=a}}
[b]: \  \ell |r\ \text{prime},  a\in \frac{\ell}{r}\Z/\Z>}. \]
For any $\sigma\in G_r=\Gal(\Q(\mu_r)/\Q)$,  if $\sigma(\zeta)=\zeta^x$, set $\sigma([a])=[xa]$. 
By this action, $U_r$ becomes  a $G_r$-module. The universal distribution is 
well known to be a free abelian group, moreover, for any integer $r'|r$, then 
the natural map from $U_{r'}$ to $U_r$ is a split monomorphism and thus $U_{r'}$ can 
be considered as  a submodule of $U_r$. 

The theory of  the universal distribution 
plays an important role in the theory of cyclotomic fields. Detailed 
information can be found in the well-known textbooks by Lang~\cite{Lang1} and 
Washington~\cite{Washington}.  Most notably,  Kubert~\cite{Kubert1}  and 
\cite{Kubert2},  and Sinnott~\cite{Sinnott1} studied the $\{\pm 1\}$-cohomology 
of $U_r$, from which, Sinnott got his famous index formula about the cyclotomic 
units and the Stickelberger elements. 

Recently, Anderson~\cite{ad2} found a brand new way to compute the 
$\{\pm 1\}$-cohomology of the universal distribution $U_r$.  He discovered a 
cochain complex which is a resolution of the universal distribution.  To study 
a certain group cohomology of $U_r$,  it is therefore translated
to study a double complex related to this group cohomology. In this paper, 
we use Anderson's resolution to construct a double complex related to the
$G_r$-cohomology of $U_r$ and study the spectral sequence of the double  
complex.   Suppose that $r$ is fixed  odd squarefree integer, 
we prove the following theorem:
\begin{theoa}[Abridged Form] The cohomology group
\[ H^{n}(G_r, U_r)=\bigoplus_{r'|r}  H^{n+n_{r'}}_{r'}(G_{r}, \Z)  \]
where $n_{r'}=$number of prime factors of $r'$ and
\[ H^{n}_{r'}(G_r, \Z):=\bigcap_{\ell|r'} 
\ker(H^{n}(G_{r}, \Z)\rightarrow H^{n}(G_{r/\ell}, \Z)). \]
\end{theoa}
We shall discuss the Unabridged Form in $\S 6$.  What's
more, for any positive integer $M$ which is a common factor of $\ell-1$ over all prime factors
$\ell$ of $r$, let $\sigma_{\ell}$ be a generator of the cyclic group $G_{\ell}$ and let
\[ D_{r'}:=\prod_{\ell |r'}\,  \sum^{\ell-2}_{k=0} k\sigma^k_{\ell}, \]
then
\begin{theob} \label{theob}
The image of the family
\[ \left \{ D_{r'}\left [\sum_{\ell | r' }\frac{1}{\ell}\right ]:  \forall\  r'|r \right \} \]
in $U_r/MU_r$ is a $\Z/M\Z$-basis for $H^{0}(G_r, U_r/MU_r)$. 
\end{theob}

Theorem B has interesting application in arithmetic. We follow the line given 
in Rubin~\cite{Rubin1}. Let
$\F=\Q(\mu_m)^{+}$ be  the maximal real subfield of $\Q(\mu_m)$,
assume  $\{\ell: \ell |r\}$ is a family of distinct odd primes which 
split completely in $\F/\Q$ and are $\equiv 1\pmod M$ for fixed integer $M$.
Suppose that we have  
a $G_r$ -homomorphism $\xi$ from  $U_r$  to $\F(\mu_r)^{\times}$.
 Then $\xi$ induces a map
\[ H^n(\xi):   H^{n}(G_r, U_r/MU_r)\longrightarrow 
H^{n}(G_r, \F(\mu_r)^{\times}/\F(\mu_r)^{\times\, M})  \]
for each $n\in \Z_{\geq 0}$.  In the case $n=0$,  since 
$H^{0}(G_r, \F(\mu_r)^{\times}/\F(\mu_r)^{\times\, M})= \F^{\times}/\F^{\times M}$, 
we have the map
\[ H^0(\xi):   H^{0}(G_r, U_r/MU_r)\longrightarrow \F^{\times}/\F^{\times M}. \]
In particular, let $\Q^{ab}$ be the abelian closure of $\Q$. 
Let $\mathbf e$ be an injective homomorphism from $\Q/\Z$ 
to  $\Q^{ab \,  \times}$.   Put
\[ \xi([a])=(\mathbf e(a+\frac{1}{m})-1)(\mathbf e(a-\frac{1}{m})-1). \]
Then $\xi$ is a $G_r$-homomorphism from $U_r$ to $\F(\mu_r)^{\times}$. 
The image  $H^0(\xi)(D_{r'}[\sum_{\ell |r'}\frac{1}{\ell}] )$ is just the Kolyvagin 
element $\kappa_{r'}$ as given in \cite{Rubin1}.  From this point of view, we can
regard the Euler system  as a system in the cohomology group   
$H^{0}(G_r, U_r/MU_r)$.  This is the initial motivation for this paper.

This paper is organized in the following order.  We give general notations
in $\S 2$.  In $\S 3$,  we study Anderson's 
resolution in detail. In $\S 4$, a special $G_r$-projective resolution 
$\bp_{\bullet}$ of $\Z$ is constructed and the group cohomology of $\Z$ and of $\Z/M\Z$ 
are given. With this projective resolution $\bp_{\bullet}$,   we construct a double complex  
in $\S 5$ whose total cohomology is 
the $G_r$-cohomology of $U_r$. The standard spectral sequence method is then
used to compute the cohomology group $H^{\ast}(G_r,U_r)$.  
In $\S 6$, we study the lifting problem and prove Theorem B. 
\vskip 0.3cm
\noindent {\textbf {Acknowledgment.} }\   I thank my advisor Professor Greg  W. Anderson for 
introducing me to this interesting field, for his insight and insistence  driving 
me to this paper.  The double complex method is the brain child of Professor 
Anderson. I also thank Dr. Hans Uli Walther for some useful discussion about Lemma~\ref{la2}.

\section{Notations} Fix a finite set $S$ of cardinality $|S|=s$.  Fix a family 
$\{\ell_i: i\in S\}$ of distinct odd prime numbers. Fix a positive integer $M$ 
dividing $\ell_i-1$ for all $i\in S$. Fix a total order  $\omega$ of $S$.  Put
\begin{itemize}
\item $r=r_S:=\prod_{i\in S} \ell_i$,
\item $G_S:=\Gal(\Q(\mu_r)/\Q)$.
\end{itemize}
For each $i\in S$, put
\begin{itemize}
\item  $G_i$:=the inertia subgroup of $G_S$ at $\ell_i$,
\item $\sigma_i$:=a fixed generator of $G_i$,
\item $N_i:=\sum_{k=0}^{\ell_i-2} \sigma_i^k$, 
$D_i:=\sum_{k=0}^{\ell_i-2} k\sigma_i^k$, 
\item $Fr_i$:=the arithmetic Frobenius automorphism at $\ell_i$ in $G_S/G_i$.
\end{itemize}
For each subset $T\subseteq S$, put
\begin{itemize}
\item $r_T:= \prod_{i\in T} \ell_i$, $\mu_T:=\mu_{r_T}$,
\item $G_T:= \prod_{i\in T} G_i\subset G_S$, 
\item $N_T:=\prod_{i\in T} N_i$, $D_T:=\prod_{i\in T} D_i$.
\end{itemize}
\noindent Put $R:=\Z_{\geq 0}[S]$.   For any element  $e=(e_i)\in R$, put
\begin{itemize}
\item $\deg e:={\sum_i} e_i$,  
\item $\supp e:=\{i\in S: e_i\neq 0\}$.
\item $\we:=(\we_i)\in R$ where $\we_i=\sum_{j<_{\omega} i} e_j$.
\end{itemize}
\noindent For any $e, e'\in R$, put $\omega(e,e'):=\sum_{j<_{\omega} i} e'_j e_i$.

For $a\in \Q/\Z$,  the order of $a$(denoted by $\ord a$) means its order in 
$\Q/\Z$. For any set $X$,  the cardinality of $X$ is denoted by $|X|$, 
the free abelian group generated by $X$ is denoted by $<\, X\, >$, 
the free $\Z/M\Z$-module generated by $X$ is denoted by $<\, X\, >_M$.  
The family of all subsets of $X$ is denoted by $2^X$. We call
a subfamily $\mi$ of $2^X$ an $order\ ideal$ of  $X$ if for all $Y\in \mi$,  
$2^Y\subseteq \mi$. For any pair of sets $X$ and $Y$, the difference 
of $X$ and $Y$ is denoted by $X\backslash Y$.  

For any complex $C^{\bullet}$, the complex $C^{\bullet}[n]$ is the complex
with components $C^{m}[n]=C^{m+n}$. For any complex $C^{\bullet}$ of
$\Z$-modules,  $C^{\bullet}_M:=C^{\bullet}\otimes \Z/M\Z$. 

\section{Universal ordinary distribution and its structure}
\subsection{Universal ordinary distribution and Anderson's resolution} 
Let $\{[a]: a\in \Q/\Z\}$ be a basis of a free
abelian group $\A$. Recall by Kubert~\cite{Kubert1}, the (rank $1$) universal 
ordinary distribution $U$ is given by 
\[ U=\frac{<[a]: a\in \Q/\Z>}{<[a]-{\displaystyle\sum_{n b=a}}
[b]:\  n\in \N>}. \]
For any positive number $f$, the universal ordinary distribution of level $f$ is given by
\[ U_f=\frac{<[a]: a\in \frac{1}{f}\Z/\Z>}{<[a]-{\displaystyle\sum_{p b=a}}
[b]: p\mid f, a\in \frac{p}{f}\Z/\Z>}. \]
For any $\sigma\in \Gal(\Q(\mu_f)/\Q)$,  set $\sigma([a])=[xa]$ if $\sigma$ sends
each $f$-th root of unity to its $x$-th power. 
By this action,  $U_f$ is a $G_f=\Gal(\Q(\mu_f)/\Q)$-module.  
Much has been studied about the structures of $U$ and $U_f$, we list some
basic properties here(for detailed proof, see Anderson~\cite{ad1, ad2},
Kubert~\cite{Kubert1, Kubert2} and Washington~\cite{Washington}). First recall
for any $a\in \frac{1}{f}\Z/\Z$, $a$ can be written uniquely
\[ a=\sum_{p\mid f}\, \sum_{\nu\in \N} \frac{a_{p\nu}}
{p^{\nu}}\pmod \Z, 0\leq a_{p\nu}\leq p-1. \]
Then 
\begin{prop} \label{fact}
1). The universal ordinary distribution $U_f$ is a free abelian group of 
rank $| G_f|$, the set $\{[a]: a\in \frac{1}{f}\Z/\Z, \ a_{p1}
\neq p-1,\  \forall p\mid f\}$ is a $\Z$-basis for $U_f$.

2). For any factor $g$ of $f$, the natural map from $U_g$ to $U_f$ is a split 
monomorphism. Moreover, by this natural map, $U$ is the direct limit of $U_f$ 
for $f\in \N$ and thus $U$ is free. 
\end{prop}
In the sequel, for our convenience, the universal distribution $U_r$ will 
be written as $U_S$ and $U_{r_T}$ as $U_T$.  Now let
 \[ \LL^{\bullet}_S=< [a, T]:  a\in \frac{r_T}{r_S}\Z/\Z, T\subseteq S>  \]
be the free abelian group generated by the symbols $[a,T]$, and let
\[ L^{p}_S=<[a, T]: |T|=-p, a\in \frac{r_T}{r_S}\Z/\Z, T\subseteq S> , \]
then $\LL^{\bullet}_S$ is a bounded graded module.  Furthermore, for any $\sigma\in G_S$, 
set $\sigma[a,T]=[xa,T]$ if $\sigma$ sends each $r$-th root of unity to its $x$-th power.
By this action,  $\LL^{\bullet}_S$ becomes a $G_S$-module.  Let 
\[ d[a, T]=\sum_{i\in T}\omega(i, T) ([a, T\backslash \{i\}]-
\sum_{\ell_i b=a}[b, T\backslash \{i\}]). \]
where
\[ \omega(i,T)=\begin{cases} (-1)^{|\{j\in T: j<_{\omega} i\}|},\ &\text{if}\
 i\in T,\\ 0, &\text{if}\ i\notin T. \end{cases}\]
It is easy to check that $d^2=0$ and $d$ is $G_S$-equivariant. Thus $(\LL^{\bullet}_S, d)$ 
is a cochain complex.  Note that the definition of $d$ depends on $\omega$. 
We'll write $d_{\omega}$ instead if we need to emphasize the order $\omega$.
The following proposition is given by Anderson:
 \begin{prop}  \label{anderson}
The $n$-th cohomology of the complex $(\LL^{\bullet}_S, d)$ is $0$ 
for $n\neq 0$ and $U_S$ for $n=0$, furthermore, the map from $L^0_S$ to
$U_S$ is given by $\mathfrak u: [a,\emptyset]\mapsto [a]$. \end{prop}
\begin{rem} 1. The above proposition(in a more general form suitable for a
resolution for
the distribution $U_f$), though  known by Anderson for quite a while, has 
no published proof by now. We put the proof in Appendix A, traces of the 
idea behind the proof can be found in Anderson~\cite{ad1, ad2}. 

2. For the sake of this proposition, we call $(\LL^{\bullet}_S,d)$ 
Anderson's resolution of the universal distribution $U_S$. This
resolution has been used by Das~\cite{Das} in his work about the
algebraic $\Gamma$-monomials and double covering of cyclotomic fields. 
\end{rem}

\subsection{Double complex structure of $\LL^{\bullet}_S$} 
A remarkable fact about
Anderson's resolution $\LL^{\bullet}_S$ is that it possesses an even more 
delicate double complex structure, which in turn gives a natural filtration
for the universal distribution $U_S$.  We start with a more careful look at
$\LL^{\bullet}_S$,  which we'll denote by $\LL^{\bullet}$ instead .  For any 
$T\subseteq S$,   we always regard  $\LL^{\bullet}_T$ as a 
subcomplex of $\LL^{\bullet}$.  Moreover, for any order ideal $\mi$ of $S$,  put
\[ \LL^{\bullet}(\mi):=\sum_{T\in \mi} \LL^{\bullet}_T\  \text{and}\  
U_S(\mi):=\sum_{T\in \mi} U_T. \]
In particular,  let $\mi(n)$ be the order ideal consisting of all subsets $T$ such that
$|T|\leq n$,  and let $\LL^{\bullet}(n)= \LL^{\bullet}(\mi(n))$ and $U_S(n)=U_S(\mi(n))$.
Note that $\LL^{\bullet}(2^T)=\LL^{\bullet}_T$.
For any $a\in \frac{1}{r_S}\Z/\Z$, let 
\[ \supp a:=\{i: \ell_i\mid \ord a\}\subseteq S.  \]
We see that 
\[ \LL^{\bullet}=<[a, T]:  \supp a\cap T=\emptyset>.  \]
Then $ \LL^{\bullet}(\mi)$ is the free abelian group generated by
\[ \left \{ [a,T]: T\cup \supp a\in \mi, \ T\cap \supp a=\emptyset \right \}, \]
and  $U_S(\mi)$ is the free abelian group generated by 
\[ \left \{[a]: \supp a\in \mi,\ a_{\ell_i 1}\neq \ell_i-1 
\ \text{for all }\ i\in S\right \}. \]
Immediately we have 
\begin{prop} \label{intersection}
Let $\mi_1$ and $\mi_2$ be two order ideals  of $S$, then

1).  $\LL^{\bullet}(\mi_1\cap \mi_2)=\LL^{\bullet}(\mi_1)\cap \LL^{\bullet}(\mi_2)$,\ 
$U_S(\mi_1\cap \mi_2)=U_S(\mi_1)\cap U_S(\mi_2)$. 

2).  $\LL^{\bullet}(\mi_1\cup \mi_2)=\LL^{\bullet}(\mi_1)+ \LL^{\bullet}(\mi_2)$,\ 
$U_S(\mi_1\cup \mi_2)=U_S(\mi_1) + U_S(\mi_2)$. 

\end{prop}
\noindent Then  
\begin{prop} \label{exact} The complex $\LL^{\bullet}(\mi)$ is acyclic
with the $0$-cohomology $U_{S} (\mi)$.
\end{prop}
\begin{proof} We let $\tL^{\bullet}( \mi)$ be the complex
\[ 0\longrightarrow \LL^{\bullet}(\mi) \stackrel{\mathfrak u}\longrightarrow U_{S}(\mi) 
\longrightarrow 0. \]
Hence it suffices to show that $\tL^{\bullet}( \mi)$ is exact.  Let $T$ be a maximal
element in the order ideal $\mi$. Let $\mi'$ be the order ideal whose maximal 
element set is the maximal element of $\mi$ excluding $T$, then 
\[ \mi=\mi'\cup 2^T. \]
By Proposition~\ref{intersection}, we have
\[ \tL^{\bullet}(\mi)/\tL^{\bullet}(2^T)=
\tL^{\bullet}( \mi')/ \tL^{\bullet}(\mi'\cap 2^T).
\]
Now we prove the Proposition 
by induction on the cardinality of maximal elements of $\mi$. If $\mi$ has only one maximal
element, this is just Proposition~\ref{anderson}.  In general, both $\mi'$ and $\mi' \cap 2^T$ have
less maximal elements than $\mi$ has. Thus the exactness of $\tL^{\bullet}(\mi)$ follows from
the exactness of the three complexes $\tL^{\bullet}(2^T)$, $\tL^{\bullet}( \mi')$ and 
$\tL^{\bullet}(\mi'\cap 2^T)$. 
\end{proof}

Now we can construct a double 
complex whose total single complex is $(\LL^{\bullet}, d)$.  With abuse of notation,
we'll write it as $\LL^{\bullet,\bullet}$. For any pair of subsets $T'$, $T$ of 
$S$ such that $T'\supseteq T$, set
\[ L_{T', T}:=<[a, T]: \supp a={S\backslash T'}>, \]
then $L_{T', T}$ is isomorphic to $\Ind^{G_S}_{G_{T'}}\Z$. Moreover, 
for any $i\in T$, the map 
\[ \varphi_i: L_{T', T}\rightarrow L_{T', T\backslash \{i\}},\  [a,T]\mapsto 
[a,T\backslash \{i\}] \]
defines a natural isomorphism between $L_{T', T}$ 
and $L_{T', T\backslash \{i\}}$.  Now for any $T\subseteq S$, 
\[ \LL^{\bullet}_{T}= \bigoplus_{T_1, T_2} L_{T_1, T_2},\ \text{where}\ 
T_2\cup (S\backslash T_1)\subseteq T,\ 
T_2\cap (S\backslash T_1)=\emptyset. \]
and if let $\Gamma(\mi):=\{(T_1,T_2): T_2\cup (S\backslash T_1)\in \mi, \ 
T_2\cap  (S\backslash T_1)=\emptyset \}$, then
\[ \LL^{\bullet}(\mi)= \bigoplus_{(T_1, T_2)\in \Gamma(\mi)} L_{T_1, T_2}. \]
In general for any $i\in S$, define
\[ \varphi_i: L^{p}\rightarrow L^{p+1},  [a,T]\mapsto \chi_T(i) 
[a,T\backslash \{i\}] \]
where $\chi_T$ is the characteristic function of $T$. 
Let $\varphi(L^{p})$ be the subgroup of $L^{p+1}$ generated by $\varphi_i(L^{p})$ 
for all $i\in S$, inductively, let $\varphi^n(L^p)$ be the subgroup of $L^{p+n}$
generated by $\varphi_i(\varphi^{n-1}(L^p))$ for all $i\in S$. By this setup,
there is a filtration of $L^{p}$ given by
\[ \varphi^{s+p}(L^{-s})\subseteq \varphi^{s+p-1}(L^{-s+1})\subseteq \cdots
\subseteq L^{p}. \]
This filtration enables us to define the double complex structure of $\LL^{\bullet}$.
For the element $[a,T]\in \LL^{\bullet}$, 
we say $[a, T]$ is of bidegree $(p_1,p_2)$ if 
$ [a,T]\in \varphi^{p_2}(L^{p_1})\backslash \varphi^{p_2+1}(L^{p_1-1})$,  more explicitly, if 
\[ p_1=|\supp a|-s, \ p_2=s-|\supp a|-|T|. \]
Then we see that all elements of $L_{T', T}$ are of bidegree $(-|T'|, |T'|-|T|)$. 
Let $L^{p_1,p_2}$ be the subgroup of $\LL^{\bullet}$ generated by all symbols
$[a,T]$ with bidegree $(p_1, p_2)$, then
\[ L^{p_1,p_2}= \bigoplus_{ |T|=-p_1-p_2}\, \bigoplus_{\substack{ |T'|=-p_1\\T'\supseteq T}}
\,  L_{T', T}.  \]
Set 
\[ d_1: L^{p_1,p_2}\rightarrow L^{p_1+1,p_2}, [a,T]\mapsto -\sum_{i\in T}
\omega(i, T) N_i[Fr^{-1}_i a +\frac{1}{\ell_i}, T\backslash \{i\}] , \]
\[ d_2:  L^{p_1,p_2}\rightarrow L^{p_1,p_2+1}, [a,T]\mapsto \sum_{i\in T}
\omega(i, T)(1-Fr^{-1}_i) [ a, T\backslash \{i\}] . \]
It is easy to check that $d^2_1=d^2_2=d_1d_2+d_2d_1=0$.  Hence we construct a 
double complex $(\LL^{\bullet,\bullet}; d_1, d_2)$. Note that $d=d_1+d_2$ and 
\[ L^p=\bigoplus_{p_1+p_2=p} L^{p_1, p_2}, \]
thus $(\LL^{\bullet},d)$ is the single total complex of the double complex 
$(\LL^{\bullet,\bullet}; d_1, d_2)$, with the second filtration  
given by $\varphi$.  Thus the total cohomology of $(\LL^{\bullet,\bullet}; d_1, d_2)$
is  the cohomology of $(\LL^{\bullet},d)$.
\begin{prop} The $E_1$ term of the spectral sequence  arising from the double 
complex $(\LL^{\bullet,\bullet}; d_1, d_2)$ 
by the first filtration(i.e., $H^{p_1}_{d_1}(\LL^{\bullet, p_2})$) is 
\[  E^{p_1, p_2}_1=\begin{cases} U_{S}(s-{p_2})/U_{S}( s-{p_2-1}),\ 
&\text{if}\  p_1=-p_2; \\ 0,\ &\text{otherwise}. \end{cases}
\]
Thus the spectral sequence for the first filtration degenerates at $E_1$. 
\end{prop}
\begin{proof} Note that
\[ \LL^{\bullet}(n)=\bigoplus_{p_2\geq s-n} L^{p_1,p_2} \]
then it is easy to see that $\LL^{\bullet, p_2}[-p_2]$ is nothing but
the quotient complex $\LL^{\bullet}(s-{p_2})/\LL^{\bullet}( s-{p_2-1})$.
The short exact sequence
\[ 0\longrightarrow  \LL^{\bullet}( s-{p_2-1})\longrightarrow
\LL^{\bullet}( s-{p_2})\longrightarrow \LL^{\bullet, p_2}[-p_2]
\longrightarrow 0 \]
induces a long exact sequence
\[ \cdots
\rightarrow H^{i}(\LL^{\bullet}(s-{p_2})) \rightarrow
H^i(\LL^{\bullet, p_2}[-p_2])\rightarrow H^{i+1}(\LL^{\bullet}( s-{p_2-1}))
\rightarrow \cdots \]
By Proposition~\ref{exact}, 
for $i\neq 0$ and $-1$, both $H^{i}(\LL^{\bullet}( s-{p_2}))$ and
$H^{i+1}(\LL^{\bullet}( s-{p_2+1}))$ are $0$, so is $H^i(\LL^{\bullet, p_2}[-p_2])$.
Therefore  the above long exact sequence is just the  exact sequence
\[ 0\rightarrow H^{-1}(\LL^{\bullet, p_2}[-p_2])\rightarrow U_{S}( s-{p_2-1})
\rightarrow U_{S}( s-{p_2}) \rightarrow H^{0}(\LL^{\bullet, p_2}[-p_2])
\rightarrow 0, \]
Since the map from $U_{S}(s-{p_2-1})$ to $U_{S}( s-{p_2})$ is injective,
the proposition follows immediately.
\end{proof}
The following results follow immediately from Proposition~\ref{fact},  Proposition~\ref{anderson}
and Proposition~\ref{exact}:
\begin{prop} 1). One has
\[ H^n(\LL^{\bullet}_M)
=\begin{cases} U_S/MU_S, & if \ n=0;\\ 0,& if \ n\neq 0. \end{cases} \]

2). Moreover, for any order ideal $\mi$ of $S$, one has
\[ H^n(\LL^{\bullet}_M(\mi))
=\begin{cases} U_S(\mi)/MU_S(\mi), & if \ n=0;\\ 0,& if \ n\neq 0. \end{cases} \]
\end{prop}

\section{The cohomology groups $H^{\ast}(G_T, \Z)$ and $H^{\ast}(G_T, \Z/M\Z)$ }
\subsection{A projective resolution of $\Z$}
We first have a convention here:  Let 
$X$ be a finite totally ordered set and $x\in X$.  Suppose $A_x$ is a module 
related to $x$.  We call 
\[ A_X=A_{x_1} \otimes \cdots \otimes A_{x_n} \]
the standard tensor product of $A_x$ over $X$ if $X=\{x_1,\cdots, x_n\}$ and 
$x_1<\cdots<x_n$. Similarly, we can define the 
standard tensor product of 
elements $a_x\in A_x$ and of complexes $A^{\bullet}_x$. 

Let 
\[ (\bp_{i\bullet}, \partial_{i}):\ \  \cdots \stackrel{
\partial_{i, j+1}}{\longrightarrow} P_{i,j+1}\stackrel{ 
\partial_{ij}}{\longrightarrow} P_{ij} \cdots \stackrel{\partial_{ 
i0}}{\longrightarrow} P_{i0}\longrightarrow  0 \]
with $P_{ij}=\Z[G_i]$ for any $j\geq 0$, $\partial_{ij}$ is the multiplication by 
$1-\sigma_{i} $ if $j$ is even and by $N_{i}$ if $j$ is odd.  It is well known 
that $\bp_{i\bullet}$ is a $\Z[G_i]$-projective resolution of trivial module $\Z$.  
For any $T\subseteq S$,  let  $\bp_{T\bullet}$ be the standard tensor product of 
$\bp_{i\bullet}$ over $i\in T$.  It is well known by homological algebra that
$\bp_{T\bullet}$ is a $\Z[G_T]$-projective resolution of trivial module $\Z$.  
Now for the collection $\{P_{i, e_i}: i\in T\}$,  the standard product of $P_{i, e_i}$ 
over $T$ is a rank $1$ free $\Z[G_T]$-module whose grade is $\sum_i e_i$.  Now let 
$e\in R$ be the element whose $i$-th component is $e_i$ if $i\in T$ and $0$ if not, 
and write the standard product of $P_{i, e_i}$ over $T$ as $\Z[G_T][e]$, then
\[ \bp_{T\bullet}=\bigoplus_{\supp e\subseteq T} \Z[G_T][e]. \]  
For any $x= (\cdots \otimes x_i \otimes \cdots) \in \Z[G_T][e]$,  the differential is 
given by
\[ \partial_T(x)=\sum_{i\in T} (-1)^{\we_i} (\cdots \otimes \partial_{i,e_i-1}
(x_i)\otimes \cdots) ,\]
In particular, for $T=S$.  let
\[ \bp_{\bullet}=\bp_{S \bullet}=\bigoplus_{e\in R} \Z[G_S] [e]. \]  
For any $T'\subseteq T$,  we have a natural inclusion $\iota: \Z[G_T'][e]
\hookrightarrow \Z[G_T][e]$
for any $e\in R$ such that $\supp e\subseteq T'$.  By this inclusion, 
$\bp_{T'\bullet}$ becomes a subcomplex of $\bp_{T\bullet}$. 

Now we define a diagonal map $\Phi_T: \bp_{T\bullet}\rightarrow  \bp_{T\bullet}\otimes  
\bp_{T\bullet}$.  First set 
\[ \begin{split} \Phi_{ie_i, ie'_i}: P_{i, e_i+e'_i}&\longrightarrow  P_{ie_i}
\otimes  P_{ie'_i}\\
1&\longmapsto \begin{cases} 1\otimes 1,\ & \text{if $e_i$ even};\\
1\otimes \sigma_i,\ & \text{if $e_i$ odd, $e'_i$ even};\\
{\displaystyle\sum_{0\leq m<n\leq \ell_i-2}}
\sigma^m_i\otimes \sigma^n_i,\ & \text{if $e_i$ odd, $e'_i$ odd};
\end{cases}
\end{split} \]
Then the map $\Phi_i: \bp_{i\bullet}\rightarrow  \bp_{i\bullet}\otimes  \bp_{i\bullet}$ 
given by $\Phi_{ie_i, ie'_i}$ is the diagonal map  for the cyclic group $G_i$
(see Cartan-Eilenberg~\cite{CE}, P250-252). For any $e,e'\in R$ with 
support contained in $T$, consider the standard product $P_{e,e'}$ of 
$P_{ie_i}\otimes P_{ie'_i}$ over $i\in T$. The isomorphism
\[ \begin{split} 
\alpha: P_{ie_i}\otimes P_{je'_j} &\longrightarrow P_{je'_j}\otimes P_{ie_i} \\
x\otimes y&\longmapsto (-1)^{e_i e'_i} y\otimes x
\end{split} \]
induces an isomorphism $\alpha: P_{e,e'}\rightarrow \Z[G_T][e]\otimes \Z[G_T][e']$ by
\[ \left (\cdots (x_i\otimes y_i)\cdots\right )\longmapsto
(-1)^{\omega(e,e')}(\cdots x_i \cdots)\otimes (\cdots y_i \cdots). \]
On the other hand, the standard product of the diagonal maps 
$\Phi_{ie_i, ie'_i}$ over $i\in T$
defines a map $\beta: \Z[G_T][e+e']\rightarrow P_{e,e'}$. We let
$\Phi_{e,e'}=\alpha\circ \beta$ and let
\[ \Phi_{T,p,q}=\sum_{\substack{e, e': \deg\ e=p, \deg\ e'=q\\ 
\supp e+e'\subseteq T}}\Phi_{e, e'}. \]
Then $\Phi_T$ defines the diagonal map from $\bp_{T\bullet}$ to
$\bp_{T\bullet}\otimes \bp_{T\bullet}$. This map enables us to 
compute the cup product structures.

\subsection{The cohomology groups $H^{\ast}(G_T,\Z)$ and $H^{\ast}(G_T,\Z/M\Z)$}
Let $\bc_i=\Hom_{G_i}(\bp_{i \bullet}, \Z)$,  then $\bc_i$ is the complex
\[  \Z\stackrel{0}{\longrightarrow }\Z
\stackrel{\ell_i-1}{\longrightarrow }\Z
\stackrel{0}{\longrightarrow }\Z \stackrel{\ell_i-1}{\longrightarrow } \cdots \]
with the initial term at degree $0$. We denote by $C^j_{i}$ the $j$-th term of 
$\bc_i$. By the theory of group cohomology, 
\[ H^{\ast}(G_i, \Z)=H^{\ast} (\bc_i). \]
Now for any $T\subseteq S$,  let
$\bc_T$ be the standard tensor product of $\bc_i$ for $i\in T$. 
If write
\[  \Hom_{G_T}(\Z[G_T][e], \Z)=\Z[e], \]
then
\[ \bc_T=\Hom_{G_T}(\bp_{T\bullet}, \Z)=\bigoplus_{\supp e\subseteq T} \Z[e]. \]  
and
\[ H^{\ast}(G_T, \Z)=H^{\ast}(\bc_T). \]
Moreover,  for any $T'\subseteq T$,  the inclusion $\iota: \bp_{T'\bullet}\hookrightarrow \bp_{T\bullet}$
induces a map
\[ \iota^{\ast}:  \bc_T \longrightarrow  \bc_{T'} , \]
which is just the  natural projection of 
\[ {\displaystyle \bigoplus_{\supp e\subseteq T}} \Z[e]\longrightarrow 
{\displaystyle \bigoplus_{\supp e\subseteq T'}} \Z[e].  \]
On the other hand, $G_{T'}$ can also be considered naturally as a quotient group
of $G_{T}$, by this meaning,  the inflation map is just the injection
\[  \bigoplus_{\supp e\subseteq T'} \Z[e] \hookrightarrow 
\bigoplus_{\supp e\subseteq T} \Z[e]. \]
Now for any $j\in \Z_{\geq 0}$ even,  let 
\[ {\mathbf C}^{\bullet j}_{i}=\begin{cases}
\cdots 0\longrightarrow C^0_i\longrightarrow 0\cdots, \ &if \ j=0;\\
\cdots 0\longrightarrow C^{j-1}_i \stackrel{\ell_i-1}\longrightarrow C^{j}_i
\longrightarrow 0\cdots, \ &if \ j>0. \end{cases} \]
For any $e=(e_i)\in 2R$, i.e., $e_i$ even for all $i\in S$,  we let $\bc_{e}$ be the standard product
${\mathbf C}^{\bullet e_i}_{i}$ over $i\in S$.  If $\supp e\subseteq T$, 
then $\bc_{e}$ is a subcomplex of $\bc_T$ and 
\[ \bc_{T}=\bigoplus_{\substack{e\in 2R\\ \supp e\subseteq T}} \bc_{e}. \]
Figure 1  shows us how the decomposition looks like in the case $S=\{1,2\}$.
Denote by $A_e$ the cohomology group $H^{\ast}(\bc_{e})$ and 
$A^n_e$ its $n$-th component.  Then 
\[ H^{\ast}(G_T,\Z)=\bigoplus_{\substack{ e\in 2R\\\supp e\subseteq T}} A_e,\ 
H^{n}(G_T,\Z)=\bigoplus_{\substack{e\in 2R\\ \supp e\subseteq T}} A^n_e.  \]

\begin{figure} 
\begin{center}
\includegraphics[angle=270]{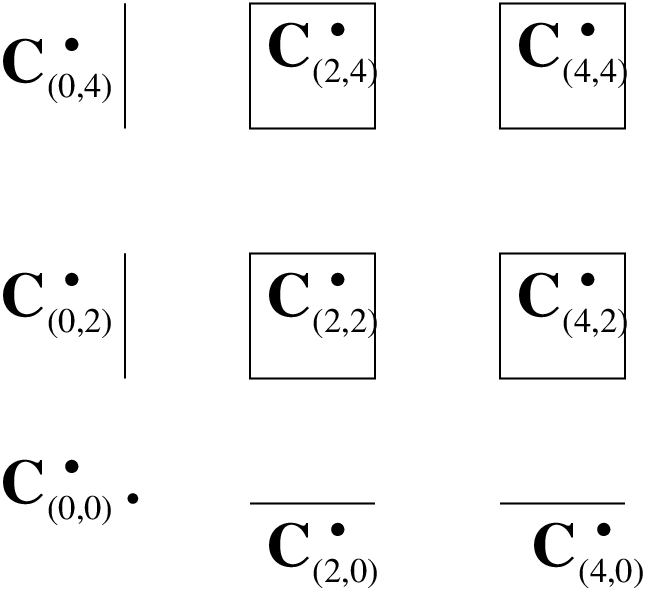}
\caption{The complex $\bc_S$ when $S=\{1,2\}$.}

\end{center}
\end{figure}

We now study the abelian group $A_e$. First we need a lemma from linear algebra:
\begin{la} \label{la1}
Let $v=(m_1, m_2, \cdots, m_n)^t$ be a $n$-dimensional column vector with integer 
entries $m_i$,  then the greatest common divisor of $m_i$ is $1$ if and only if 
there exists an $n\times n$ matrix $A\in SL_n(\Z)$ whose first column is $v$.
\end{la}

Now suppose $\supp e=T=\{i_1, \cdots, i_t\}$ and $|T|=t$. 
If $t=0$, then $T=\emptyset$,  it is easy to see that $A_e=A^0_e=\Z$. 
Now if $T\neq \emptyset$,
we claim that $\bc_e[\deg e-t]$ is isomorphic to the exterior algebra $\Lambda(x_1,
\cdots, x_t)$ with differential $d(x)=\sum (\ell_i-1) x_i\wedge x$ and $\deg x_i=1$. 
This claim is easy to check:  First if $t=1$,  let $T=\{i\}$,  then ${\mathbf C}^{\bullet e_i}_{i}=
C^{e_i-1}\oplus C^{e_i}$.  This case is trivial. 
In general,  if  ${\mathbf C}^{\bullet e_i}_{i}[e_i-1]$ is isomorphic to $\Lambda(x_i)$, 
the tensor product of  ${\mathbf C}^{\bullet e_i}_{i}[e_i-1]$ is nothing but 
$\bc_e[\deg e-t]$ and
the tensor product of $\Lambda(x_i)$ is just $\Lambda(x_1,
\cdots, x_t)$, hence they are isomorphic to each other.

Now let $m_T$ be the greatest common divisor of $\ell_i-1$ for $i\in T$, 
thus the greatest common divisor of $(\ell_i-1)/m_T$ is 1, let $A$ be the matrix 
given by Lemma \ref{la1} corresponding to the vector $(\cdots,.(\ell_i-1)/m_T, 
\cdots)$.  Let $(y_1,\cdots , y_t)=(x_1,\cdots, x_t)A$. Then $\{y_1, \cdots, 
y_t\}$ is a set of new generators for the above exterior algebra and we have 
$d(x)=m_T\, y_1\wedge x$.  We see easily that
\[ H^{\ast}(\Lambda(x_1,\cdots, x_t))=(\Z/m_T \Z )^{2^{t-1}} \]
and 
\[ H^{j}(\Lambda(x_1,\cdots, x_t))=(\Z/m_T \Z )^{\binom{t-1}{j}},  0\leq j\leq t-1. \]
Combining the above analysis, we have
\begin{prop}  \label{propa} 
There exists a family of complexes
\[ \{ \bc_e \subseteq  \bc=\Hom_{G_S}(\bp_{\bullet}, \Z): e\in 2R\} \]
such that 

1). For each $T\subseteq S$, we can identify $\bc_T=\Hom_{G_T}(\bp_{T \bullet}, \Z)$ with
${\displaystyle  \bigoplus_{\substack{ e\in 2R\\ \supp e \subseteq T}} }\bc_e$
through the following splitting exact sequence
\[ 0 \longrightarrow \bigoplus_{\substack{ e\in 2R\\ \supp e \not \subseteq T}} \bc_e
\longrightarrow \bc  \longrightarrow \bc_T \longrightarrow 0 \]

2).  The cohomology groups $H^{\ast}(\bc_e)=A_e$ and
$H^{n}(\bc_e)=A^n_e$ are given by:
\begin{quote} 
(a).  If $\supp e\neq \emptyset$,  let $m_e$ be the greatest common divisor of $\ell_i-1$ 
for $i\in \supp e$,    then $A_e$ 
is the abelian group $(\Z/m_e \Z )^{2^{|\supp e|-1}}$,  and 
\[ A^n_e=\begin{cases} (\Z/m_e 
\Z )^{\binom{|\supp e|-1}{j}}, \ & \text{ if} \  n=\deg e -j\ \text{and} \  0\leq j\leq |\supp e|-1;\\  
0,\ & \text{if otherwise}. \end{cases} \]  
(b).  If $\supp e=\emptyset$, then $A_e=A^0_e=\Z$. 
\end{quote}
\end{prop}

For the case $H^{\ast}(G, \Z/M\Z)$,  the situation is much easier.  We have
\begin{prop} \label{propb}
There exists a family 
\[ \left \{ [e]\in H^{\ast}(G_S, \Z/M\Z): e\in R\right \} \]
with the following properties:

1). For each $T\subseteq S$ and $n\in \Z_{\geq 0}$, the restriction of the 
family 
\[ \left \{[e]: e\in R, \ \supp e\subseteq T,\ \deg e=n\right \} \]
to $H^n(G_T, \Z/M\Z)$ is a $\Z/M\Z$-basis of the latter.

2). For each $T\subseteq S$ and $e\in R$ such that $\supp e\nsubseteq T$, 
the restriction of $[e]$ to $H^{\ast}(G_T, \Z/M\Z)$ vanishes.

3). One has the cup product structure in $H^{\ast}(G_T, \Z/M\Z)$ given by
\[ [e]\cup [e']=(-1)^{\omega(e,e')} \prod_{\substack{i\in S\\e_ie'_i\equiv 1(2)}} 
\left ( \frac{\ell_i-1}{2}\right ) [e+e'] \]
for all $e, e'\in R$.

\end{prop}

\begin{proof}
 The complex ${\mathbf C}^{\bullet}_{M,i}=\Hom_{G_i}(\bp_{i \bullet}, \Z/M\Z)$
by definition,  is a complex 
with ${\mathbf C}^{j}_{M,i}=\Z/M\Z$ for $j\geq 0$ and the differential $0$.
In general,  ${\mathbf C}^{\bullet}_{M,T}=\bc_T\otimes \Z/M\Z$ is exactly
 the standard tensor product of  ${C}^{\bullet}_{M,i}$ for all $i\in T$.  Write 
 \[  {\mathbf C}^{\bullet}_{M,T}=\Hom_{G_i}(\bp_{T\bullet}, \Z/M\Z)
=\sum_{\supp e \subseteq T} \Z/M\Z[e]. \]
Since now ${\mathbf C}^{\bullet}_{M,T}$ has differential $0$,  $H^{\ast}({\mathbf C}^{\bullet}_{M,T})
={\mathbf C}^{\bullet}_{M,T}$.  The restriction map is easy to see. This finishes the proof
of 1) and 2).

For the cup product, the diagonal map $\Phi_T$  
given above naturally induces a map:
\[ {\mathbf C}^{\bullet}_{M,T}\times {\mathbf C}^{\bullet}_{M,T}
\longrightarrow   {\mathbf C}^{\bullet}_{M,T} \]
which defines the cup product structure. More specifically, the cup product
map
\[ \Z/M\Z[e] \times \Z/M\Z[e'] \longrightarrow \Z/M\Z[e+e'] \]
is induced from $\Phi_{e,e'}$.  Now the claim follows soon from the 
explicit expression of $\Phi_{e,e'}$.
\end{proof}

\section{Study of $H^{\ast}(G_S, U_S)$}
\subsection{The complex $\K$}  With the preparation in $\S 3$ and $\S 4$, 
set
\[  \K^{\bullet,\bullet}:=\Hom_{G_S}(\bp_{\bullet}, \LL^{\bullet}).  \] 
Let $d$ and $\delta$ be the induced differentials of $d$ and $\partial$ by 
$\LL^{\bullet}$ and $ \bp_{\bullet}$ respectively.  If we let
\[ [a,T,e]:=(e\mapsto [a,T])\in \Hom_{G_S}(P_{e}, <[a,T]>), \]
then 
\[ \begin{split} 
&K^{p,q} =<[a,T,e]:  a\in \frac{r_T}{r_S}\Z/\Z,  |T|=-p,   \deg e=q>;\\
&d[a,T,e]=\sum_{i\in T} \omega(i, T)([a,T\backslash \{i\}, e]- \sum_{\ell_i b=a}
[b, T\backslash \{i\}, e]); \\
&\delta[a,T,e]=(-1)^{|T|} \sum_{i\in S} (-1)^{\we_i}\cdot \begin{cases} 
(1-\sigma_i)[a,T, e+\epsilon_i], \ &if \ e_i\ even;\\N_i 
[a,T, e+\epsilon_i], \ &if\ e_i \ odd. \end{cases} \end{split} \]
For any $T\subseteq S$, set
\[ \K^{\bullet,\bullet}(T)=\Hom_{G_S}(\bp_{\bullet}, \LL^{\bullet}_T)=
<[a, T', e]: [a, T']\in \LL^{\bullet}_T,   e\in R>\]
and 
\[ \K^{\bullet,\bullet}_T=\Hom_{G_T}(\bp_{T \bullet}, \LL^{\bullet}_T)
=<[a, T', e]: [a, T']\in \LL^{\bullet}_T,   e\in R, \supp e\subseteq  T >. \] 
Furthermore, for any order ideal $\mi$, set
\[ \K^{\bullet,\bullet}(\mi):=\Hom_{G_S}(\bp_{\bullet}, \LL^{\bullet}(\mi))=\sum_{T\in \mi} 
\K^{\bullet,\bullet}(T). \]
and set 
\[ \K^{\bullet,\bullet}(n):=\Hom_{G_S}(\bp_{\bullet}, \LL^{\bullet}(n)).\]
Set
\[ \fu^{\bullet}:=\Hom_{G_S}(\bp_{\bullet}, U_S)=\frac{<[a, e]: a\in \frac{1}{r_S}\Z/\Z, e\in R>}
{<[a,e]-\sum_{\ell_i b=a} [b,e]: a\in \frac{\ell_i}{r_S}\Z/\Z, e\in R>}, \]
with the induced differential $\delta$ from $\partial$.  Correspondingly, 
\[ \fu^{\bullet}(\mi):=\frac{<[a, e]: a\in \frac{1}{r_T}\Z/\Z\ \text{for some}\ T\in \mi, e\in R>}
{<[a,e]-\sum_{\ell_i b=a} [b,e]: a\in \frac{\ell_i}{r_T}\Z/\Z\ \text{for some}\ T\in \mi, e\in R>}, \]
which is a subcomplex of $\fu^{\bullet}$.  We consider $\fu^{\bullet}$
as the double complex
$(\fu^{\bullet, \bullet}; 0, \delta)$ concentrated on the vertical axis. From 
Proposition~\ref{anderson}, we have a map
\[ \mathfrak u: \K^{\bullet,\bullet}\rightarrow \fu^{\bullet, \bullet},\ 
[a,T,e]\mapsto \begin{cases} [a,e], \ & \text{if}\ T=\emptyset; \\
 0,\  & \text{if}\ T\neq \emptyset. \end{cases} \]

\begin{prop} \label{xin}
The map $\mathfrak u$(resp. its restriction) is a quasi-isomorphism between 
$\K^{\bullet,\bullet}$ (resp. $\K^{\bullet,\bullet}(\mi)$)
and $\fu^{\bullet, \bullet}$(resp. $\fu^{\bullet,\bullet}(\mi)$). Therefore

1). $ H^{\ast}_{total}(\K^{\bullet,\bullet})=H^{\ast}(G_S, U_S)$, 
 $H^{\ast}_{total}(\K^{\bullet,\bullet}(\mi))=H^{\ast}(G_S, U_S(\mi))$.

2). $ H^{\ast}_{total}(\K^{\bullet,\bullet}_M)=H^{\ast}(G_S, U_S/MU_S)$, 
 $H^{\ast}_{total}(\K^{\bullet,\bullet}_M(\mi))=H^{\ast}(G_S, U_S(\mi)/MU_S(\mi))$.

\end{prop}
\begin{proof} Immediately from Proposition~\ref{anderson}(resp. Proposition~\ref{exact}
for $\mi$),   we see that $\ker \mathfrak u$
is $d$-acyclic, by spectral sequence argument, it is thus $(d+\delta)$-acyclic.  On the
other hand, $\mathfrak u$ is surjective. Thus $\mathfrak u$ is a quasi-isomorphism.
Now 1) follows directly from the quasi-isomorphism.  For 2), just consider $\mathfrak u\otimes 1$,
which is also a  quasi-isomorphism.
\end{proof}

From Proposition~\ref{xin}, the $G_S$-cohomology of $U_S$  is isomorphic to
the total cohomology of the double complex $(\K^{\bullet,\bullet}; d, \delta)$
Therefore we can  use the spectral sequence of the double complex $\K^{\bullet,\bullet}$
to study  the $G_S$-cohomology of $U_S$.  The spectral sequence of $\K^{\bullet,\bullet}$
from the second filtration exactly gives us Proposition~\ref{xin}.  Now we study the 
spectral sequence from the first filtration. Then
\[ {E}^{p,q}_1(\K^{\bullet,\bullet})=H^q_{\delta}({ \K}^{p,\bullet})
=H^q(G_S, L^p). \]
Now since 
\[  L^p=\bigoplus_{p_1+p_2=p} L^{p_1,p_2}=
\bigoplus_{|T|=-p}\, \bigoplus_{T'\supseteq T} {L}_{T',T}, \]
then 
\[ E^{p,q}_1(\K^{\bullet,\bullet})=\bigoplus_{|T|=-p}\, \, 
\bigoplus_{T'\supseteq T}H^q(G_S, L_{T', T}). \]
For the double complex $\K^{\bullet,\bullet}(\mi)$, Recall
$\Gamma(\mi)=\{(T_1,T_2): T_2\cup (S\backslash T_1)\in \mi, \ 
T_2\cap  (S\backslash T_1)=\emptyset \}$ in $\S3$, we have
\[ E^{p,q}_1(\K^{\bullet,\bullet}(\mi))=\bigoplus_{(T',T)
\in \Gamma(\mi)} H^q(G_S, L_{T', T}). \]

\subsection{A Lemma } Suppose that  for any $T\subseteq S$, there is an abelian group 
$B_T$ related to $T$,  set
 \[ A_T=\bigoplus_{T''\subseteq T} B_{T''}. \]
Then for any $T'\supseteq T$, there is a natural projection from $A_{T'}$ 
to $A_T$. Now  let $\cc^{\bullet}_{S,T}$  be  the cochain complex with components given by
\[ \cc^n_{S,T}=
\bigoplus_{\substack{|T'|=s-n\\ T'\supseteq (S\backslash T)}} A_{T'},  \]
and differential $d$ given by
\[ \begin{split} d: \  A_{T'}&\longrightarrow \bigoplus_{i\in T'\cap T} 
A_{T'\backslash \{i\}}\\
x&\longmapsto \sum_{i\in T'\cap T} \omega(i, T'\cap T) x|_{T'\backslash \{i\}}, \end{split} \]
where $x|_{T'\backslash \{i\}}$ is the projection of $x$ in $A_{T'\backslash \{i\}}$.
It is easy to verify that $\cc^{\bullet}_{S,T}$ is indeed a chain complex.  Furthermore, we 
have
\begin{la} \label{la2} For any $T\subseteq S$, 
\[  H^{n}(\cc^{\bullet}_{S,T},d)=\begin{cases} \bigoplus_{T'\supseteq T} B_{T'}, &\text{if}\  
n=0;\\ 0,  &\text{otherwise}. \end{cases} \]
\end{la}
\begin{proof}  Let $\mc^{\bullet}_{S,T}$ be the subcomplex of $\cc^{\bullet}_{S,T}$
 with the same components as 
$\cc^{\bullet}_{S,T}$ except at degree $0$,  where 
\[ \mc^0_{S,T}=\bigoplus_{T'\not \supseteq T} B_{T'}. \]
We only need to show that $\mc^{\bullet}_{S,T}$ is  exact.
We show it by double induction to the cardinalities of $S$ and $T$. 
If $T=\emptyset$,  we get a trivial complex. If $S$ consists of only one element,  
or if $T$ consists only one element, it is also trivial to verify.  In general,  
suppose  $i_0=\max\{i:i\in T\}$.  Let $S_0=S\backslash \{i_0\}$ and 
$T_0=T\backslash \{i_0\}$.
Then we have the following commutative diagram which is exact on the columns:
\[ \begin{CD}
0 @>>> \mc^0_{S,T_0}  @>{\bar{d}}>> \mc^1_{S,T_0} @>{\bar{d}}>>\cdots
\mc^{t-1}_{S,T_0}@>>>0@>>>0\\
  @.  @AA{p}A  @AA{p}A @AA{p }A@AA{p }A@. \\
0 @>>> \mc^0_{S,T}  @>{d}>> \mc^1_{S,T} @>{d}>>\cdots
\mc^{t-1}_{S,T}@>d>>\mc^{t}_{S,T}@>>>0\\
 @.  @AA{i}A  @AA{i}A @AA{i }A @AA{i }A @. \\
0 @>>> {\underset{S_0\supseteq T'\supseteq T_0} \bigoplus}B_{T'}   @>{d}>> A_{S_0} 
@>{d}>>\cdots \mc^{t-1}_{S_0,T_0}@>d>>\mc^{t}_{S_0,T_0}@>>>0
\end{CD}  \]
Here $p$ means projection and $i$ means inclusion.  The differential $\bar{d}$ is 
induced by the differential $d$  of the second row.  Notice that the third row 
is a variation of the chain complex $\mc^{\bullet}_{S_0,T_0}$, the first row is the chain 
complex $\mc^{\bullet}_{S,T_0}$. By induction, the first row and and the third row are 
exact, so is the middle one. 
\end{proof}
\subsection{The Study of $E_2$ terms} 
By $\S 4.1$, we know that 
\[ E^{p,q}_1(\K^{\bullet,\bullet})=\bigoplus_{|T|=-p}\, \, 
\bigoplus_{T'\supseteq T} H^q(G_S, L_{T', T}). \] 
Now let's consider the induced differential $\bar d$ of $d$ in the $E_1$ term. 
Since $d=d_1+d_2$, we can also write $\bar d = \bar d_1+\bar d_2$.  We first 
look $\bar d_2$, which is induced by
\[ \begin{split}
{L}_{T',T}&\longrightarrow \bigoplus_{i\in T} { L}_{T',T\backslash \{i\}}\\
[a,T] &\longmapsto \sum_{i\in T}\omega(i,T) (1-Fr^{-1}_i)[a,T\backslash \{i\}]
\end{split} \]
Since for any $i\in T$, ${L}_{T',T}$ and ${ L}_{T',T\backslash \{i\}}$ are 
$G_S$-isomorphic by the map $\varphi_i$, and since for any $q\geq 0$, $H^q(G, A)$ 
is a trivial $G$-module, we have
\[  \bar d_2=\sum_{i\in T}\omega(i,T)(1-Fr^{-1}_i) \bar{\varphi_i}=0. \]
For the map $\bar d_1$, which is induced by
\[ \begin{split}
 {L}_{T',T}&\longrightarrow  \bigoplus_{i\in T} {L}_{T'\backslash \{i\},
T\backslash \{i\}} \\
[a,T] &\longmapsto  - \sum_{i\in T}\omega(i,T) N_i[Fr^{-1}_i a+
\frac{1}{\ell_i}, T\backslash \{i\}]  
\end{split} \]
For any $i\in T$,  let 
\[ \begin{split} \psi_i: {L}_{T',T}&\longrightarrow {L}_{T'\backslash \{i\},
T\backslash \{i\}}\\ [a, T]&\longmapsto N_i[Fr^{-1}_i a+\frac{1}{\ell_i}, 
T\backslash \{i\}]  \end{split} \]
$\psi_i$ is a $G_S$-homomorphism and therefore induces a map in the $G_S$-cohomology:
\[ H^q(\psi_i): H^q(G_S, {L}_{T',T}) \rightarrow H^q(G_S, 
{L}_{T'\backslash \{i\},T\backslash \{i\}}). \]
We have the commutative diagram:
\[ \begin{CD}
 L_{T',T}@>{\psi_i}>> {L}_{T'\backslash \{i\}, T\backslash \{i\}}\\
 @VV{\theta_{T'}}V @VV{\theta_{T'\backslash \{i\}} }V\\
\Z@>{res}>>\Z
\end{CD} \]
where the top row are $G_S$-modules,  the left $\Z$ is a trivial 
$G_{T'}$-module and the right $\Z$ a trivial $G_{T'\backslash \{i\}}$-module,
and $\theta_{T'}$ is the homomorphism 
sending $[\frac{1}{r_{S\backslash T'}}, T]$ to $1$ and $[\frac{x}{r_{S\backslash T'}}, T]$
to $0$ if $x\neq 1$. 
Then the above diagram induces the following commutative diagram:
\[ \begin{CD}
 H^q(G_{S},  L_{T', T})@>{H^{q}(\psi_i)}>>
H^q(G_{S},  {L}_{T'\backslash \{i\}, T\backslash \{i\}}) \\
 @VV{\theta^{\ast}_{T'}}V @VV{\theta^{\ast}_{T'\backslash \{i\}} }V\\
 H^q(G_{T'}, \Z)@>{res}>> H^q(G_{T'\backslash \{i\}}, \Z)
\end{CD} \]
where $\theta^{\ast}_{T'}$(and $\theta^{\ast}_{T'\backslash \{i\}} $) is the isomorphism 
given by Shapiro's lemma(See 
Serre~\cite{Serre2}, Chap. VII, $\S 5$, Exercise). 
We identify $H^q(G_{S},  L_{T', T})$ with $H^q(G_{T'}, \Z)$,
moreover,   to keep track on $T$, we'll write $H^q(G_{T'},  \Z)$ 
as $H^q(G_{T', T}, \Z)$. Then we see that $H^q(\psi_i)$ is the restriction map 
from $H^q(G_{T', T}, \Z)$
to $H^q(G_{T'\backslash \{i\}, T\backslash \{i\}}, \Z)$.  
The induced differential $\bar d=\bar d_1$ is exactly the map
\[ \begin{split} H^q(G_{T',T}, \Z)&\longrightarrow \bigoplus_{i\in T} 
H^q(G_{T'\backslash \{i\}, T\backslash \{i\}}, \Z)\\
x&\longmapsto -\sum_{i\in T} \omega(i,T) x_i \end{split} \]
where $x_i$ is the restriction of $x$ in $H^q(G_{T'\backslash \{i\}, 
T\backslash \{i\}}, \Z)$. Hence  we have a cochain complex $\cc(q;S, T)$
\[ \  H^q(G_{S,T},\Z){\stackrel{\bar d_1}\longrightarrow} 
\bigoplus_{i\in T} H^q
(G_{S\backslash \{i\}, T\backslash \{i\}},\Z)\cdots
{\stackrel{\bar d_1}\longrightarrow }
H^q(G_{S\backslash T,\emptyset},\Z)\longrightarrow 0\]
Note that the complex $E^{\bullet, q}_1(\K^{\bullet, \bullet})$ is just the direct sum of $\cc(q;S, T)$
over all subsets $T$ of $S$. Moreover, the complex $E^{\bullet, q}_1(\K^{\bullet, \bullet})(\mi)$
is the direct sum of $\mc(q;S, T)$ over all subsets $T\in {\mi}$.

Recall in Proposition~\ref{propa}, we obtained
\[ H^q(G_{T},\Z)=\bigoplus_{\substack{e\in 2R\\ \supp e\subseteq T}} A^q_e. \]
If let
\[ A^q_T=H^q(G_{T},\Z), \  B^q_T=\bigoplus_{\substack{e\in 2R\\ \substack{\supp e=T}}} A^q_e. \]
then we have $A^q_T=\bigoplus_{T''\subseteq T} B^q_{T''}$.  The complex
$\cc(q; S,T)[-|T|\, ]$satisfies the conditions in Lemma \ref{la2},  
thus the  $n$-th cohomology of the cochain complex $\cc(q; S,T)$
is $0$ if $n\neq -|T|$ 
and $\underset{T'\supseteq T}{\sum}B_{T'}$ if $n=-|T|$. We have the following proposition:
\begin{prop}  \label{e2}
One has

1). $  E^{p,q}_2(\K^{\bullet,\bullet})\cong\displaystyle{ \bigoplus_{|T|=-p}\, \, 
\bigoplus_{\substack{e\in 2R \\ \supp e\supseteq T}}} A^q_e $.   

2). 
$  E^{p,q}_2(\K^{\bullet,\bullet}(\mi))\cong \displaystyle{ \bigoplus_{\substack{|T|=-p\\T\in {\mi}}}\, \, 
\bigoplus_{\substack{e\in 2R\\ \supp e\supseteq T}}} A^q_e $.  
\end{prop}

\subsection{Proof of Theorem A}  Finally we are in a stage to 
prove Theorem A.  Put 
\[ \fS^{\bullet,\bullet} =<[a,T,e]\in \K^{\bullet,\bullet}, a\neq 0\ \text{if} \ 
\supp e\supseteq T>. \]
It is easy to verify that $\fS^{\bullet,\bullet}$ is a subcomplex of 
$\K^{\bullet,\bullet}$. Set 
\[ \fQ^{\bullet,\bullet} =\K^{\bullet, \bullet}/\fS^{\bullet,\bullet}=<[0,T,e]: 
\supp e\supseteq T>. \]
Note that the induced differential of $d$ at $\fQ^{\bullet,\bullet}$ is $0$.  Moreover, set 
\[ \fS^{\bullet,\bullet}(\mi):=\K^{\bullet,\bullet}(\mi)\cap \fS^{\bullet,\bullet}, \] 
and 
\[  \fQ^{\bullet,\bullet}(\mi):=\K^{\bullet,\bullet}(\mi)/\fS^{\bullet,\bullet}(\mi)=<[0, T,e]:\ 
T\in \mi, \  \supp e \supseteq T>. \]
Let $f$ be the corresponding quotient map, then we  have a commutative diagram:
\[ \begin{CD}
 \K^{\bullet, \bullet}_M(\mi)@>{inc}>> \K^{\bullet, \bullet}_M\\
 @VV{f}V @VV{f }V\\
\fQ^{\bullet, \bullet}_M(\mi)@>{inc}>> \fQ^{\bullet, \bullet}_M
\end{CD} \]
We claim that
\begin{prop} \label{theoq}
The quotient map $f: \K^{\bullet,\bullet}\rightarrow \fQ^{\bullet,\bullet}$ is a 
quasi-isomorphism. Moreover, 
the quotient map $f: \K^{\bullet,\bullet}(\mi)\rightarrow \fQ^{\bullet,\bullet}(\mi)$
is a quasi-isomorphism.
\end{prop}
\begin{proof}  Let  
\[ \hl_T:=<[0, T, e]: e\in R>=\Hom_{G_S}(\bp_{\bullet}, L_{S,T})\subseteq \K^{\bullet, \bullet} \]
and let
\[ \al_T:=<[0, T, e]:  \supp e\supseteq T>,\   \bl_T:=<[0, T, e]:  \supp e\subseteq S\backslash\{i\},\
\text{for some}\ i\in T>\]
Through the map $ L_{S,T}\rightarrow \Z, \ [0,T]\mapsto 1$, we have a commutative diagram
\[ \begin{CD}
\hl_T  @{=}  \al_T@{.} {\displaystyle \bigoplus} @{.} \bl_T\\
    @VVV  @VVV @{.} @VVV\\
\bc@{=} {\displaystyle \bigoplus_{\substack{e\in 2R\\\supp e\supseteq T}}}  \bc_e\ @{.} 
{\displaystyle \bigoplus} @{.} 
\displaystyle{\bigoplus_{\substack{i\in T, \ e\in 2R\\ \supp e\subseteq  S\backslash \{i\}}}}  \bc_e
\end{CD}  \]
where $\bc$ and $\bc_e$ are given in Proposition~\ref{propa}.  By this diagram, we
identify $\hl_T$ with $\bc$.  By Proposition~\ref{propa}, we have
\[ \ker( H^{\ast}(G_S, \Z)\rightarrow  H^{\ast}(G_{S\backslash \{i\}}, \Z) )
=H^{\ast}( \bigoplus_{\substack{e\in 2R\\ i\in \supp e}} \bc_e). \]
Then by the proof of Proposition~\ref{e2},
\[ \begin{split}  \ker(\bar d|_{H^{q}(\hl_T)})=& \bigcap_{i\in T} 
\ker( H^{\ast}(G_S, L_{S, T})\rightarrow  H^{\ast}(G_{S}, L_{S\backslash \{i\}, T\backslash \{i\}}) )\\
=&  \bigcap_{i\in T} H^{\ast}( \bigoplus_{\substack{e\in 2R\\ i\in \supp e}} \bc_e)
= H^{\ast}( \bigcap_{i\in T} \bigoplus_{\substack{e\in 2R\\ i\in \supp e}} \bc_e )\\
=& H^{\ast}(\bigoplus_{\substack{e\in 2R\\ T\subseteq \supp e}} \bc_e )
= H^{\ast}(\al_T) \end{split} \]
where the second and the last identities we use the isomorphisms given in the above diagram.
Hence we have
\[ E^{p,q}_2(\K^{\bullet, \bullet})=\bigoplus_{|T|=-p}\ker(\bar d|_{H^{q}(\hl_T)})=
\bigoplus_{|T|=-p} H^{q}(\al_T). \]
On the other hand,  
\[  \fQ^{\bullet,\bullet}=\bigoplus_{T\subseteq S}  \al_T. \]
Since $d=0$ in $\fQ^{\bullet,\bullet}$, the spectral sequence of $\fQ^{\bullet,\bullet}$ by the
first filtration(i.e., by $d$) degenerates at $E_1$. We have
\[  E_1^{p,q}( \fQ^{\bullet,\bullet})=E_2^{p,q}( \fQ^{\bullet,\bullet})
=\bigoplus_{|T|=-p} H^{q}(\al_T). \]

Since the projective map from $\hl_T$ to $\al_T$ in the commutative diagram is nothing
but the restriction of the quotient map $f$ at $\hl_T$, by the above analysis,  we
get an isomorphism
\[ f_2: E^{p,q}_2(\K^{\bullet, \bullet})\longrightarrow E_2^{p,q}( \fQ^{\bullet,\bullet}). \]
Thus the spectral sequence of $\K^{\bullet, \bullet}$ and $\fQ^{\bullet,\bullet}$ are
isomorphic at $E_r$ for $r\geq 2$. In our case, the first filtration is finite, thus
strongly convergent,  therefore $f$ is a quasi-isomorphism(see 
Cartan-Eilenberg~\cite{CE}, Page 322, Theorem 3.2). 

The case $\mi$ is similar. In this case, 
\[ E^{p,q}_2(\K^{\bullet, \bullet})(\mi)=\bigoplus_{\substack{T\in \mi\\|T|=-p}}
\ker(\bar d|_{H^{q}(\hl_T)})=
\bigoplus_{|T|=-p} H^{q}(\al_T), \]
and 
\[  \fQ^{\bullet,\bullet}(\mi)=\bigoplus_{T\in \mi}  \al_T. \]
Now follow the same analysis as above.
\end{proof}
For any subset $T$ of $S$,  set
\[ H^{\ast}_T(G_S, \Z):=\bigcap_{i\in T}\ker(H^{\ast}(G_S, \Z)\rightarrow 
H^{\ast}(G_{S\backslash \{i\}}, \Z)) \]
we see that 
\[ H^{\ast}(\al_T)\cong H^{\ast}_T(G_S, \Z) \]
by the identification of $\hl_T$ and $\bc$.  The following theorem is the main result
in the paper:
\begin{theoa}[Unabridged Form] \label{theoa} 
1).  The cohomology group $H^{\ast}(G_S, U_S)$
is given by
\[ H^{\ast}(G_S, U_S)=\bigoplus_{T\subseteq S} H^{\ast}_T(G_S, \Z)[\, |T|\, ]
=\bigoplus_{T\subseteq S} \bigoplus_{\substack{e\in 2R \\ \supp e\supseteq T}} A_e[\, |T| \, ]. \]
where $A_e[\, |T| \, ]$ represents the cohomology group $H^{\ast}(\bc_e[\, |T| \, ])$.

2).  The cohomology group $H^{\ast}(G_S, U_S(\mi))$
is given by
\[ H^{\ast}(G_S, U_S(\mi))=\bigoplus_{T\in \mi} H^{\ast}_T(G_S, \Z)[\, |T| \, ]
=\bigoplus_{T\in \mi} \bigoplus_{\substack{e\in 2R\\ \supp e\supseteq T}} A_e[\, |T| \, ]. \]
\end{theoa}
\begin{proof}  We only prove 1).  The proof of 2) follows the same route.
By Proposition~\ref{xin} and Proposition~\ref{theoq}, we know that
\[ H^{\ast}(G_S, U_S)=H^{\ast}_{total}(\K^{\bullet})=H^{\ast}_{total}(\fQ^{\bullet}). \]
Now 
\[ H^{n}_{total}(\fQ^{\bullet})= \bigoplus_{T\subseteq S} H^{n+|T|}(\al_T). \]
Part 1) follows immediately.
\end{proof}
\begin{rem} We can see that Part 1) is actually a special case of Part 2)  when the order
ideal $\mi$ is $2^S$.
\end{rem}
In the case $\Z/M\Z$,  we have
\begin{theo} \label{theom}
There exists a family
\[ \left \{ c_{\text{\tiny T},e}\in H^{\ast}(G_S, U_S/MU_S): T\subseteq S,\ e\in R, \ 
\supp e\supseteq T\right \} \]
with the following properties:

1). For each $n\in \Z_{\geq 0}$, the subfamily
\[ \left \{ c_{\text{\tiny T},e}: T\subseteq S,\ e\in R,\ \supp e\supseteq T, \deg e=n+|T|\right \} \]
is a $\Z/M\Z$-basis for $H^{n}(G_S, U_S/MU_S)$.

2). For any order ideal $\mi$ of $S$,  let $U_S(\mi)=\sum_{T\in \mi} U_T$.
By the inclusion $U_S(\mi)\hookrightarrow U_S$, $H^{\ast}(G_S, 
U_S(\mi)/ MU_S(\mi))$ can be considered as a submodule of $H^{\ast}(G_S, U_S/M_S)$.
Furthermore, the subfamily
\[ \left \{ c_{\text{\tiny T},e}: T\in {\mi},\ e\in R,\ \supp e\supseteq T\right \} \]
is a $Z/MZ$ basis for $H^{\ast}(G_S, U_S(\mi)/ MU_S(\mi))$.

3). One has cup product structure
\[ [e']\cup  c_{\text{\tiny T},e}=(-1)^{\omega(e', e)}\prod_{\substack{i\in S\\e_ie'_i\equiv
1(2)}} \left ( \frac{\ell_i-1}{2}\right ) c_{\text{\tiny T}, e+e'} \]
for all $e,e'\in R$ and $T\subseteq \supp e$.
\end{theo}
\begin{proof}
1).  By Proposition~\ref{theoq}, we have induced quasi-isomorphism:
\[ f\otimes 1: \K^{\bullet, \bullet}_M\longrightarrow 
\fQ^{\bullet, \bullet}_M  \]
Now since the induced differentials of $d$ and $\delta$ in 
$\fQ^{\bullet, \bullet}_M$ are $0$. 
Consider the cocycle $[0, T, e]$
in $\fQ^{\bullet,\bullet}_M$, there exists a cocycle 
$C_{T,e}$(unique modulo boundary) which is the lifting of $[0, T, e]$ 
by the quotient map $f\otimes 1$.  Hence 
 $\mathfrak u(C_{T,e})\otimes 1$ is a cocycle in
the complex $\fu^{\bullet}_M$. Let  $c_{\text{\tiny T},e}$ denote the cohomology 
element in $H^{\ast}(G_S, U_S/MU_S)$ represented by the cocycle  $\mathfrak u(C_{T,e})\otimes 1$.
Then $\{ c_{\text{\tiny T},e}: \supp e\supseteq T\}$ 
is a canonical $\Z/M\Z$-basis for the cohomology group $H^{\ast}(G_S, U_S/MU_S)$. 
This finishes the proof of 1).

2). Similar to 1),  just consider the map $f\otimes 1: \K^{\bullet, \bullet}_M(\mi)\rightarrow 
\fQ^{\bullet, \bullet}_M(\mi)$.

3). For the cup product,  there is natural homomorphism
\[ \Z/M\Z \otimes U_S/MU_S \longrightarrow  U_S/MU_S, \]
therefore $H^{\ast}(G_S, U_S/MU_S)$(and also $H^{\ast}(G_S, U_S(\mi)/MU_S(\mi))$
has a natural $H^{\ast}(G_S, \Z/M\Z)$-module 
structure. By the theory of spectral sequences(see, for example 
Brown~\cite{Brown}, Chap. $7$,  $\S 5$),  we have the cochain cup product 
\[ \bc_M \otimes \K^{\bullet, \bullet}_M \longrightarrow 
\K^{\bullet, \bullet}_M. \]
By using the diagonal map $\Phi_S$ defined in $\S 3$, it is easy to check that:
\[ \bc_M \otimes \fS^{\bullet,\bullet}_M
\subseteq \fS^{\bullet,\bullet}_M, \]
hence we can pass the cup product structure to the quotient and have
\[  \bc_M \otimes \fQ^{\bullet,\bullet}_M
\longrightarrow \fQ^{\bullet,\bullet}_M. \]
Now 3) follows immediately from the explicit expression 
of $\Phi_S$. This concludes the proof.
\end{proof}

\section{Explicit basis of  $H^{0}(G_S, U_S/MU_S)$}
In $\S 5$, we obtain a canonical basis $\{c_{\text{\tiny T},e}: \supp \supseteq T\}$ for the 
cohomology group $H^{\ast}(G_S, U_S/MU_S)$.  However,  little is known yet for 
the explicit expression of the cocycle $c_{\text{\tiny T},e}$ in the complex 
$\Hom_{G_S}(\bp_{\bullet}, U_S/MU_S)$, which makes it necessary to study how to 
lift the cocycle $[0, T, e]$ in $\fQ^{\bullet,\bullet}_{M}$ to the cocycle $C_{T,e}$ 
in $\K^{\bullet,\bullet}_M$. Unfortunately,  we are unable to get a complete answer
for this problem in this paper.  We obtain partial solution in the $0$-cocycles case,
however,  which is enough for us to prove Theorem B. 

\subsection{The triple complex structure of $\K$}
Recall in $\S 3$, $\LL$ has a double complex structure, therefore
we can make $\K$ as a triple complex.  Set
\[ K^{p_1,p_2,q}:=\Hom_{G_S}(\bp_{\bullet}, L^{p_1, p_2})=<[a,T,e]: [a,T]\in L^{p_1, p_2}, \deg e=q> \]
and the differentials $(d_1, d_2, \delta)$ given by
\[ d_1[a, T,e]=-\sum_{i\in T} \omega(i, T) N_i[Fr^{-1}_i a+\frac{1}{\ell_i}, 
T\backslash \{i\}, e] \]
\[ d_2[a, T,e]=\sum_{i\in T} \omega(i, T) (1-Fr^{-1}_i)[a, T\backslash \{i\}, e] \]
and $\delta$ as given in the double complex $\K^{\bullet,\bullet}$.
In this setup, we see that $\K(\mi)$ becomes a triple subcomplex of $\K$, moreover  
\[ \K(n)=\bigoplus_{p_2\geq s-n} K^{p_1, p_2, q}. \] 
Correspondingly, we define triple complex structures on $\K_M$, $\K_M(\mi)$ and
$\K_M(n)$.  This triple complex structure enables us to construct different double
complex structures in $\K$ and $\K_M$. By studying those double complexes,  
we can gather more information about $\K$. This method
will be illustrated on next subsection.

\subsection{The double complex $(\K^{\bullet, p_2, \bullet}_M, d_1, \delta)$} 
For fixed $p_2$,  let
\[ \K^{\bullet, p_2, \bullet}_M=\bigoplus_{p_1, q} K^{p_1, p_2, q}_M, \]
with differentials $d_1$ and $\delta$, then we get a double complex 
$(\K^{\bullet, p_2, \bullet}_M; d_1, \delta)$.   Similarly, we can get the 
double complex
$(\K^{\bullet, \bullet}_M; d_1+\delta, d_2)$ whose $(p_1+q, p_2)$-component is 
$\bigoplus K^{p_1,p_2, q}_M$.  As before, for any $\mi$, we have double 
complexes $\K^{\bullet, p_2, \bullet}_M(\mi)$ and $\bigoplus K^{p_1,p_2, q}_M(\mi)$
which are subcomplexes of $\K^{\bullet, p_2, \bullet}_M$ and $\bigoplus 
K^{p_1,p_2, q}_M$ respectively. First we have
\begin{prop} \label{p1q}
1). $H^{\ast}_{total}(\K^{\bullet, p_2, \bullet}_M; 
d_1, \delta)$ is a free $\Z/M\Z$-module generated by 
cocycles $C'_{T,e}$ with leading term $[0, T, e]$ and the remainder
with $q$-degree less than $\deg e$ over all   
pairs $(T,e)$ satisfying $|T|=s-p_2$ and $\supp e\supseteq T$.

2). Moreover, $H^{\ast}_{total}(\K^{\bullet, p_2, \bullet}_M(\mi); 
d_1, \delta)$ is a free $\Z/M\Z$-module generated by 
cocycles $C'_{T,e}$ with leading term $[0, T, e]$ and the remainder
with $q$-degree less than $\deg e$ over all   
pairs $(T,e)$ satisfying $T\in {\mi}$,  $|T|=s-p_2$ and $\supp e\supseteq T$.
\end{prop}
\begin{proof}  We only prove (1).  The proof of (2) is similar. 
First look  the spectral sequence of $\K^{\bullet, p_2, \bullet}_M$
with the second filtration(i.e., the filtration  given by $q$),  then 
\[ E^{p_1, q}_1(\K^{\bullet, p_2, \bullet}_M)=H^q(G_S, L^{p_1, p_2}). \]
Next for the differential $d_1$ induced on $E_1$, with the same analysis as in 
computing the
$E_2$ terms of $(\K; d, \delta)$(see $\S 4$, Proposition~\ref{e2}), we have
\[ E^{p_1, q}_2(\K^{\bullet, p_2, \bullet}_M)=\begin{cases} 
{\displaystyle \bigoplus_{|T|=s-p_2}}\, \,
{\displaystyle\bigoplus_{\substack{e:\, \deg e =q \\ \supp e\supseteq T} }}\Z/M\Z, 
\ &\text{if}\ p_1=-s;\\0,\ & \text{if}\ p_1\neq -s.  \end{cases} \]
Furthermore,    let $(\fQ^{\bullet, p_2, \bullet}_M; 0,0)$ be the double complex 
generated by all symbols $[0,T,e]$ satisfying $|T|=s-p_2$ and $ \supp e\supseteq T$, 
which can be considered as
a quotient complex  of $\K^{\bullet, p_2, \bullet}_M$.  Similar to the proof of
Theorem A,
the quotient map induces an isomorphism between the cohomology groups 
of each other.  Now let $C'_{T,e}$ be the canonical lifting of the cocycle
$[0,T,e]$ in $\fQ^{\bullet, p_2, \bullet}_M$,  then $C'_{T,e}$ is a cocycle
in $\K^{\bullet, p_2, \bullet}_M$ with the leading term $[0,T,e]$ and the remainder
contained in the direct sum of $K^{p'_1,p_2,q'}$ where $q'<\deg e$ and 
$p'_1+q'=\deg e-s$. 
\end{proof}

\begin{prop} \label{p1qp2}
The spectral sequence of the double complex $(\K^{\bullet,  
\bullet}_M; d_1+\delta, d_2)$ with the first filtration, degenerates at 
$E_1$.  The spectral sequence of the double complex $(\K^{\bullet,  
\bullet}_M(\mi); d_1+\delta, d_2)$ with the first filtration, degenerates at 
$E_1$. 
\end{prop}
\begin{proof} We only prove the first part.
The $E_1$-terms of the spectral sequence are  
\[ E^{p_1+q, p_2}_1(\K^{\bullet,\bullet}_M)
=H^{p_1+q}_{total}(\K^{\bullet, p_2, \bullet}_M; d_1, \delta). \]
Note that $|E^{p,q}_1|\geq |E^{p,q}_2|\geq \cdots\geq |E^{p,q}_{\infty}|$ in 
general for any spectral sequence, then
\[ 
|\bigoplus_{p_1+p_2+q=n} H^{p_1+q}_{total}(\K^{\bullet, p_2, \bullet}_M; d_1, \delta) | 
\geq |H^n_{total}(\K^{\bullet,  \bullet}_M, d+\delta)|. 
 \]
By Theorem A and Proposition~\ref{p1q},  the left hand side and 
the right hand side of the above inequality have the same number of elements.
hence the inequality is actually an identity.  Therefore,  the spectral sequence 
of $\K^{\bullet,\bullet}_M$ 
with filtration given by $p_1+q$ degenerates at $E_1$. 
\end{proof}

The advantage of studying the triple complex structure of the complex $\K_M$ 
is that we can obtain the $(-p_2)$-cocycles of $\K^{\bullet, p_2, 
\bullet}_M$ rather quickly.  Recall that
\[ (1-\sigma_i)D_i= N_i\pmod M. \]
Now for  the $(-p_2)$-cocycles $C'_{T,e}$,  the pair $(T, e)$ must satisfy 
$\deg e=|T|$ and therefore $e=e_T:={\displaystyle \sum_{i\in T} }\epsilon_i$. 
In this case,  for any $i\in T$,  we always have 
\[ \omega(i,T)=(-1)^{\omega(e_T)_{i}}=(-1)^{\omega(e_{T\backslash \{i\}})_{i}}. \]
First
\[ \delta [0,T,e]=0,  \ d_1 [0,T,e_T]=-\sum_{i\in T} \omega(i,T) N_i 
[\frac{r_{T\backslash \{i\}}} {\ell_i}, T\backslash \{i\}, e_T], \]
then 
\[ \delta (\sum_{i\in T} D_i [\frac{r_{T\backslash \{i\}}} {\ell_i}, T\backslash \{i\}, 
e_T\backslash \{i\}])= (-1)^{|T|}  d_1 [0,T,e_T], \]
Continue  this procedure, we have
\[  C'_{T,e_T}= \sum_{T'\subseteq T} (-1)^{|T'|(2|T|-|T'|-1)/2} 
D_{T'}[\sum_{i\in T'} 
\frac{r_{T\backslash T'}}{\ell_i}, T\backslash T', e_{T\backslash T'}]. \]

Apparently, we see that if $T\in {\mi}$, then the cocycles $C'_{T,e_T}$
 are all contained in the subcomplex $\K^{\bullet,p_2,\bullet}_M(\mi)$.
Combining the above results, we have
\begin{prop}  \label{lifting}
1). The canonical basis $\{C'_{T,e_T}: |T|=s-p_2\}$ of 
$H^{(-p_2)}(\K^{\bullet,p_2,\bullet}_M)$ is given by
\[  C'_{T,e_T}= \sum_{T'\subseteq T} (-1)^{|T'|(2|T|-|T'|-1)/2} 
D_{T'}[\sum_{i\in T'} 
\frac{r_{T\backslash T'}}{\ell_i}, T\backslash T', e_{T\backslash T'}]. \]

2). If we restrict our attention in the subcomplex $\K^{\bullet,p_2,\bullet}_M(\mi)$, 
then $H^{(-p_2)}(\K^{\bullet,p_2,\bullet}_M(\mi))$ has a canonical basis
$\{C'_{T,e_T}: |T|=s-p_2, T\in {\mi}\}$.  
\end{prop}

\subsection{Proof of Theorem B}  First we claim that 
\[ D_T[\sum_{i\in T} \frac{1}{\ell_i}]\in H^0(G_S, U_S/MU_S)=(U_S/MU_S)^{G_S}. \]
We prove it by induction on $|T|$.  For $T=\{j\}$, it is easy to see that
$(1-\sigma_i)D_j[\frac{1}{\ell_j}]=0$ for all $i\in S$. Now in general, for any $j\in T$,
\[ (1-\sigma_j) D_T[\sum_{i\in T} \frac{1}{\ell_i}]=(Fr_j-1)
D_{T\backslash \{j\}} [\sum_{i\in T\backslash \{j\}} \frac{1}{\ell_i}] \]
which is $0$ by induction, for $j\notin T$, it is obviously $0$.  Hence the claim holds.

Now we consider the double complex 
$(\K^{\bullet,\bullet}_M, d_1+\delta, d_2)$. By Proposition~\ref{p1qp2}, we know 
that $(\K^{\bullet,\bullet}_M, d_1+\delta, d_2)$ degenerates at $E_1$ for the first
filtration. By Proposition~\ref{lifting},  
$E^{-p_2,p_2}_1(\K^{\bullet,\bullet}_M)$ is generated by 
$\{C'_{T,e_T}: |T|=s-p_2\}$.  We plan to lift $C'_{T,e_T}$ to a $0$-cocycle in
$\K^{\bullet,\bullet}_M$, which is guaranteed by the degeneration at $E_1$. Moreover,
we can study the lifting $C'_{T,e_T}$ in $\K^{\bullet,\bullet}_M(T)$. 
Therefore there 
exists a cocycle $\tc_{T,e_T}$ in $\K^{\bullet,\bullet}_M(T)$ with the 
leading term $C'_{T,e_T}$ and  the 
remainder contained in the direct product of $K^{p'_1, p'_2, q'}_M(T)$ 
where $p'_1+p'_2+q'=0$ and $p'_2>p_2$.  Hence the image 
$\mathfrak u(\tc_{T,e_T})$ is exactly of the form
\[ \pm\,  D_T[\sum_{i\in T} \frac{1}{\ell_i}]+ Re(T), \]
the remainder $Re(T)=\sum n_a [a]$ where $\ord(a)$ is a proper factor 
of $r_T$.  Both $\mathfrak u(\tc_{T,e_T})$ 
and $D_T[{\displaystyle \sum_{i\in T}} \frac{1}{\ell_i}]$  are $0$-cocycles 
of $U_S/MU_S$,  so is $Re(T)$. 

In order to prove Theorem B,  it is sufficient to prove 
\[  (\ast):  Re(T)=\text{linear combination of}\ D_{T'}[\sum_{i\in T'} \frac{1}{\ell_i}] \]
for $T'$ proper subsets of $T$.  We show $(\ast)$ by induction to $|T|$. If $|T|=1$,
this is trivial. 
Now in general,  without loss of generality,  we suppose
that $T=S$ and for any $T'\subsetneq S$, $Re(T')$ is a linear combination
of $D_{T''}[{\displaystyle \sum_{i\in T''}} \frac{1}{\ell_i}]$ for $T''\subsetneqq
T'$.  Then $\mathfrak u(\tc_{T',e_T'})$  for any $ T'\subsetneq S\}$ is a 
linear combination of  $D_{T''}[{\displaystyle \sum_{i\in T''}} \frac{1}{\ell_i}]$
with $T''\subseteq T$.  By Proposition~\ref{xin},
Proposition~\ref{p1qp2} and Theorem A,  $H^0(G_S, U_S(s-1)/MU_S(s-1))$ is 
generated by  $\{\mathfrak u(\tc_{T',e_T'}): T'\subsetneq S\}$ and hence by
$D_{T'}[{\displaystyle \sum_{i\in T'}} \frac{1}{\ell_i}]$.  But obviously $Re(S)\in 
U_S(s-1)/MU_S(s-1)$, so $(\ast)$ holds for $Re(s)$. Theorem B is proved.
\begin{rem} One natural question to ask is if the bases of 
$H^{0}(G_S, U_S/MU_S)$ obtained in Theorem~\ref{theom} and in Theorem B are
the same. Unfortunately, they are not the same even in the case 
$|S|=3$. Right now, we don't know too much about the explicit 
expression of the cocycles $c_{\text{\tiny T},e}$. A deep understanding of those
cocycles should tell us more about the arithmetic of the cyclotomic 
fields.
\end{rem}

\appendix
\section{  A resolution of the universal ordinary
distribution}
{\centerline{Greg W.\ Anderson}}
\centerline{gwanders@math.umn.edu}
\vskip 0.3cm
\subsection{Basic definitions}
\subsubsection{The universal ordinary
distribution}
Let
$\AAA$ be a free abelian group equipped with a
basis $\{[x]\}$ indexed by $x\in \Q\cap [0,1)$.
For all $x\in
\Q$ put $[x]:=[\langle x\rangle]$,
where $\langle x\rangle$ is the unique rational
number in the interval $[0,1)$ congruent to $x$
modulo $1$. The {\em universal ordinary
distribution}
$U$ is defined to be the quotient
of $\AAA$ by the subgroup generated
by all elements of the form
\[[x]-\sum_{i=1}^f\left[\frac{x+i}{f}\right]
\;\;\;\left(f\in
\Z_{>0},\;\;
 x\in
\Q\right).\] 
\subsubsection{The universal ordinary
distribution of level $f$}
Fix a positive integer $f$. Let
$\AAA(f)$ be the subgroup of
$\AAA$ generated by the set
$\left\{[x]\left|x\in\frac{1}{f}\Z\right.\right\}$.
The {\em universal ordinary distribution} $U(f)$
of {\em level $f$} is defined to be the quotient
of
$\AAA(f)$ by the subgroup generated by all
elements of the form
\[[x]-\sum_{i=1}^g\left[\frac{x+i}{g}\right]
\;\;\;\left(g\in \Z_{>0}, \;g\mid f,\;\;x\in
\frac{g}{f}\Z\right).\]
The inclusions $\AAA(f)\subset \AAA$ induce a
natural isomorphism
$\lim_{\rightarrow}U(f)\iso U$.

\subsubsection{The ring $\Lambda$ and its action
on $\AAA$} Let $\Lambda$
be the polynomial ring over
$\Z$ generated by a family $\{X_p\}$ of
independent variables indexed by primes $p$, and
for each positive integer $f$, put
\[X_f:=\prod X_{p_i}^{e_i}\in \Lambda,\;\;\;
Y_f:=\prod(1-X_{p_i})^{e_i}\in \Lambda\]
where $f=\prod_i p_i^{e_i}$ is the prime
factorization of $f$. Each of the families
$\{X_f\}$ and $\{Y_f\}$ is a basis for
$\Lambda$ as a free abelian group.
We equip
$\AAA$ with
$\Lambda$-module structure by the rule
\[X_p[x]=\sum_{i=1}^p\left[\frac{x+i}{p}\right]\]
for all primes $p$ and $x\in \Q$. 
One has
\[U=\AAA/\left(\sum_p Y_p\AAA\right).\]
This last observation is the motivation for all
the results to follow.
\subsection{The structure of $\AAA$ as a
$\Lambda$-module}
\subsubsection{Partial fraction expansions}
Each $x\in \Q$ has
a unique partial fraction expansion
\[x=x_0+\sum_p\sum_{i}\frac{x_{pi}}{p^i}\]
where $p$ ranges over primes,
$i$ ranges over positive integers,
$x_0\in \Z$, $x_{pi}\in \Z\cap [0,p)$,
and all but finitely many of the coefficients
$x_{pi}$ vanish. For each
nonnegative integer $n$, put
\[\rep_n:=\left\{x\in \Q
\left|\begin{array}{l} \mbox{There exist
at most $n$ primes}\\
\mbox{$p$ such that
$x_{p1}=p-1$.}
\end{array}\right.\right\}\cap [0,1)\]
and let $\AAA_n$ be the subgroup
of $\AAA$ generated by $\{[x]\mid x\in
\rep_n\}$. 
\begin{Theorem}\label{Theorem:Free}
The following hold:
\begin{enumerate}
\item For all positive integers $f$ and $n$, one
has
\[\AAA_n\cap\AAA(f)\subseteq
\AAA_{n-1}\cap \AAA(f)+\sum_{p\mid f}
\AAA_n\cap \AAA(f/p),\]
where the sum is extended over primes $p$
dividing $f$. 
\item For each positive integer $f$, the family 
$\left\{X_g[x]\right\}$
indexed by the set
\[\left\{(g,x)\in \Z_{>0}\times \rep_0
\left|
g\mid f,\;\;x\in
\frac{g}{f}\Z\right.\right\}\] is a basis for
$\AAA(f)$.
\item For each positive integer $f$, the family 
$\left\{Y_g[x]\right\}$
indexed by the set
\[\left\{(g,x)\in \Z_{>0}\times \rep_0
\left|
g\mid f,\;\;x\in
\frac{g}{f}\Z\right.\right\}\] is a basis for
$\AAA(f)$.
\item The family $\left\{X_f[x]\right\}$
indexed by pairs $(f,x)\in \Z_{>0}\times \rep_0$
is a basis for
$\AAA$.
\item The family $\left\{Y_f[x]\right\}$
indexed by pairs $(f,x)\in \Z_{>0}\times \rep_0$
is a basis for
$\AAA$.
\item The free abelian group $\AAA$ is
free as a $\Lambda$-module,
and the family
$\left\{[x]\right\}$ indexed by $x\in \rep_0$ 
is a $\Lambda$-basis for $\AAA$.
\end{enumerate}
\end{Theorem}
\begin{proof} 1. For each $x\in \frac{1}{f}\Z\cap
(\rep_n\setminus \rep_{n-1})$, there exists
some prime $p$ dividing $f$ such that
$x_{p1}=p-1$, and one has
\[ [x]=-
\left(\sum_{i=1}^{p-1}\left[x+\frac{i}{p}\right]\right)
+X_p[px], \]
whence the result.

2. The family
$\{X_g[x]\}$ generates
$\AAA(f)$ by what we have already proved.
The family $\{X_g[x]\}$ is of
cardinality 
\[ \sum_{g\mid
f}\left|\rep_0\cap\frac{g}{f}\Z\right|
=\sum_{g\mid
f}\left|(\Z/(f/g)\Z)^\times\right|=
\sum_{g\mid
f}\left|(\Z/g\Z)^\times\right|=f. \]
Therefore the family
$\{X_g[x]\}$ is a basis for $\AAA(f)$.

3-6. These assertions follow trivially from what
we have already proved.
\end{proof}
\begin{Corollary}
The following hold:
\begin{enumerate}
\item For each $f\in \Z_{>0}$, the group
$U(f)$ is free abelian and the family
$\{[x]\}$ indexed by $x\in \frac{1}{f}\Z\cap
\rep_0$ gives rise to a basis for $U(f)$.
\item The group $U$ is free abelian and the
family
$\{[x]\}$ indexed by
$x\in
\rep_0$ gives rise to a basis for $U$.
\item The natural map $U(f)\rightarrow U$
is a split monomorphism.
\end{enumerate}
\end{Corollary}
\begin{proof} Clear. \end{proof}
\subsection{Construction of resolutions}
\subsubsection{The complex $(L,d)$}
Let $L$ be a free abelian group equipped
with a basis $\{[x,g]\}$
indexed by pairs $(x,g)$ with $x\in \Q\cap
[0,1)$ and $g$ a squarefree positive integer.
For all $x\in \Q$ and squarefree integers $g$,
put $[x,g]:=[\langle x\rangle,g]$. For all
$x\in \Q$ and increasing sequences
$p_1<\cdots<p_m$ of primes, we declare the symbol
$[x,p_1\cdots p_m]$ to be of degree $-m$ and we
set
\[\begin{array}{rcl}
&&d[x,p_1\cdots p_m]\\\\
&:=&\displaystyle\sum_{i=1}^{m}
(-1)^{i-1}
\left([x,p_1\cdots p_{i-1}p_{i+1}\cdots p_m]
-\sum_{j=1}^{p_i}
\left[\frac{x+j}{p_i},p_1\cdots
p_{i-1}p_{i+1}\cdots p_m\right]
\right),
\end{array}\]
thereby equipping the group $L$ with a grading
and a differential $d$ of degree $1$.
The map $[x,1]\mapsto [x]$ induces an isomorphism
$H^0(L,d)\iso U$.
\subsubsection{The subcomplexes $(L(f),d)$}
Fix a positive integer $f$. We define
$L(f)$ to be the graded subgroup
spanned by the symbols of the form $[x,g]$
where $g$ divides $f$ and $x\in \frac{g}{f}\Z$.
It is clear that $L(f)$ is $d$-stable.
The map $[x,1]\mapsto [x]$ induces an isomorphism
$H^0(L(f),d)\iso U(f)$.
\subsubsection{The noncommutative ring
$\tilde{\Lambda}$} Let $\tilde{\Lambda}$ be the
exterior algebra over $\Lambda$ generated by a
family of symbols
$\{\Xi_p\}$ indexed by primes $p$.
For each increasing sequence $p_1<\cdots<p_m$ of
prime numbers, put
\[\Xi_{p_1\cdots p_m}:=\Xi_{p_1}\wedge
\cdots \wedge \Xi_{p_m}\in \tilde{\Lambda},
\]
and declare $\Xi_{p_1\cdots p_m}$
to be of degree $-m$, 
thereby defining a $\Lambda$-basis
$\{\Xi_h\}$ for $\tilde{\Lambda}$ indexed by
squarefree positive integers $h$ and equipping
$\tilde{\Lambda}$ with a $\Lambda$-linear
grading. Let
$d$ be the unique
$\Lambda$-linear derivation of $\tilde{\Lambda}$
of degree $1$ such that
\[d\Xi_p=Y_p\] for all primes $p$. One then
has
\[d\Xi_{p_1\cdots p_m}=
\sum_{i=1}^{m}(-1)^{i-1} Y_{p_i}\Xi_{p_1\cdots
p_{i-1}p_{i+1}\cdots p_m}\]
for all increasing sequences $p_1<\cdots<p_m$
of prime numbers. 
\subsubsection{The subcomplexes
$(\tilde{\Lambda}(f),d)$}
Fix a positive
integer
$f$. The graded subgroup $\tilde{\Lambda}(f)$ 
generated by all elements of the
form
$Y_g\Xi_h$ where $gh$ divides $f$ is
$d$-stable. It is not difficult to verify that
the complex
$(\tilde{\Lambda}(f),d)$ is acyclic in nonzero
degree,  and  that
$H^0(\tilde{\Lambda}(f),d)$ is a free abelian
group of rank $1$ generated by the symbol
$\Xi_1=1$. 
\subsubsection{The action of $\tilde{\Lambda}$
on $L$}
We equip $L$ with graded left
$\tilde{\Lambda}$-module structure by the rules
\[\Xi_p[x,p_1\cdots
p_m]=\left\{\begin{array}{ll}
(-1)^{|\{i\mid p_i<p\}|}[x,pp_1\cdots
p_m]&\mbox{if
$p\not\in \{p_1,\dots,p_m\}$}\\
0&\mbox{if $p\in \{p_1,\dots,p_m\}$}
\end{array}\right.
\]
and
\[X_p[x,p_1\cdots p_m]=\sum_{i=1}^p
\left[\frac{x+i}{p},p_1\cdots p_m\right]\]
for all primes $p$ and increasing sequences
$p_1<\cdots<p_m$ of primes. 
By a straightforward calculation that we omit,
one can verify that
\[d(\xi\eta)=(d\xi)\eta+(-1)^{\deg \xi}\xi(d\eta)\]
for all homogeneous $\xi\in \tilde{\Lambda}$
and $\eta\in L$. 
\begin{Theorem}The following hold:
\begin{enumerate}
\item For each positive integer $f$, the complex
$(L(f),d)$ is acyclic in nonzero degree.

\item The complex $(L,d)$ is acyclic in
nonzero degree.
\end{enumerate}
\end{Theorem}
\begin{proof}  We have only to prove the first
statement.   By Theorem~\ref{Theorem:Free}
and a straightforward calculation that we omit,
one has
\[L(f)=\bigoplus_{(x,g)}\tilde{\Lambda}(g)[x,1]\]
where the direct sum is indexed by pairs
$(x,g)$ with $x\in
\frac{1}{f}\Z\cap\rep_0$ 
and $g$ is the largest positive integer such that
$x\in \frac{g}{f}\Z$. Each of the subcomplexes
$(\tilde{\Lambda}(g)[x,1],d)$ is an isomorphic
copy of $(\tilde{\Lambda}(g),d)$, and the latter
we have already observed to be  acyclic in
nonzero degree. \end{proof}
\subsubsection{Note on references}
The construction of $(L(f),d)$ presented here is
cobbled together
from ideas presented in the author's papers
\cite{ad1}
and \cite{ad2}, along with
simplifications suggested by many conversations
with Pinaki Das and Yi Ouyang on these topics.

\end{document}